\documentclass{article}

\usepackage{Jared-MK}
\usepackage{graphicx} 
\usepackage{mathtools} 
\usepackage{hyperref}
\usepackage{thmtools}
\usepackage{thm-restate}
\usepackage{cleveref}

\newcommand{\RV}{\mathcal{V}}
\newcommand{\bnu}{\overline{\nu}}

\newcommand{\bg}{\overline{g}}

\newcommand{\bh}{\overline{h}}

\newcommand{\bGamma}{\overline{\Gamma}}
\newcommand{\eps}{\epsilon}

\newcommand{\FPz}{\text{FP}_{z=0}}
\newcommand{\FPd}{\text{FP}_{\delta \to 0^+}}

\newcommand{\p}{\partial}
\newcommand{\n}{\nabla}

\newcommand{\bz}{\overline{z}}
\newcommand{\by}{\overline{z}}
\newcommand{\bB}{\overline{B}}
\newcommand{\tx}{\tilde{x}}
\newcommand{\tz}{\tilde{z}}

\newcommand{\RA}{\mathcal{A}}

\newcommand{\ty}{\tilde{y}}
\newcommand{\bJac}{\overline{\text{Jac}}}
\newcommand{\va}{v_{\alpha}}
\newcommand{\vb}{v_{\beta}}
\newcommand{\F}[2]{{}_{#1}F_{#2}}
\newcommand{\tY}{\tilde{Y}}

\newtheorem*{theorem*}{Theorem}
\newtheorem{example}{Example}

\newtheorem*{corollary*}{Corollary}

\title{An Inverse Problem for Renormalized Area: Determining the Bulk Metric with Minimal Surfaces} 

\author{Jared Marx-Kuo}
\begin{document}

\maketitle

\tableofcontents

\begin{abstract}
	\noindent We present an inverse problem connecting asymptotically hyperbolic, conformally compact metrics which are partially even to the renormalized area functional on minimal submanifolds. We use a rigidity argument to determine the conformal infinity of the metric via the renormalized area. We then consider renormalized volume of perturbations of the hemisphere to determine the higher order terms in the asymptotic expansion of the metric. We prove rigidity when these metrics are analytic, and further note that renormalized area determines the obstruction tensor for PE metrics.
\end{abstract}
\section{Introduction}
We call $(M^{n+1}, g)$ a conformally compact, asymptotically hyperbolic manifold, if the metric expands as 
\begin{align} \label{AHExpansion}
g &= \frac{dx^2 + \omega(y,x)}{x^2} \\ \nonumber
\omega(y,x) &= \omega_0(y) + \omega_1(y) x + \dots
\end{align}
near the boundary. Here, $x$ is a boundary defining function for $M$ vanishing on $\partial M$ and $\omega(y,x) \in \text{Sym}^2(T \partial M)$. Given $Y^2 \subseteq M^{n+1}$ with $\partial Y = \overline{Y} \cap \partial M$, the divergent form of the metric equation \eqref{AHExpansion} means that the area is a priori infinite. However, if we assume that $\omega_1 = 0$ in equation \eqref{AHExpansion}, then for any $\eps > 0$, one can expand
\[
\int_{x > \epsilon} dA_{Y} = a_0\epsilon^{-1} + a_2 + O(\epsilon)
\]
and define the \textbf{renormalized volume}
\[
\mathcal{V}(Y^2):= a_2
\]
Though $\mathcal{V}(Y)$ no longer represents the ``volume" of $Y$, it is a Riemannian invariant that reflects the topology and conformal geometry of $Y$ when $m$ is even (cf \cite{alexakis2010renormalized}, Proposition 3.1). When $m$ is odd, the definition depends on the choice of representative of the conformal infinity of $g$, but there is a computable ``conformal anomaly" that is of physical interest. \nl 
\indent In this work, we motivate and explore the inverse question: given the renormalized volume on minimal surfaces (or generally, submanifolds), can we recover the expansion of $g$ in equation \eqref{AHExpansion}?
%
%
\subsection{Motivation}
Renormalized volume was originally studied in high energy physics and string theory. For a $k$-brane, one can associate a $k$-dimensional submanifold, $Y$, of an ambient manifold, $\overline{X}^{n+1}$. The expected value of the Wilson line operator of the boundary, $W(\partial Y)$, is then given by $\exp(- T \RV(Y))$ where $T$ is the string tension and $\RV(Y)$ is the renormalized volume \cite{graham1999conformal}. Henningson and Skenderis \cite{henningson1998holographic} were the first to compute renormalized volume (in the literature, ``Weyl Anomaly") for low dimension odd examples, and Graham and Witten developed the mathematical theory shortly after. \nl 
\indent In their seminal work, Alexakis-Mazzeo developed a Gauss-Bonnet type formula for the renormalized area of minimal surfaces
\begin{proposition}[Alexakis-Mazzeo, Prop 3.1] \label{GaussBonnetProp}
Let $(M^{n+1}, g)$ and $\gamma \subseteq \partial M$ a $C^{3,\alpha}$ embedded curve. Suppose that $Y^2 \subseteq M^{n+1}$ is a properly embedded minimal surface with asymptotic boundary $\gamma$, then
\begin{equation} \label{GaussBonnet}
\mathcal{A}(Y) = - 2\pi \chi(Y) - \frac{1}{2} \int_Y |\hat{A}|^2 dA_Y + \int_Y \tr_Y(W_M) dA_Y
\end{equation}
where $\hat{A}$ denotes the trace-free second fundamental form and $W_M$ is the Weyl curvature tensor of $M$.
\end{proposition}
\noindent Equation \eqref{GaussBonnet} displays the conformal and topological summands which contribute to renormalized volume. An important corollary of theorem \ref{GaussBonnetProp} is the following rigidity in hyperbolic space:
\begin{corollary} \label{Rigidity}
	Let $\gamma \subseteq \R^n$ a $C^{3,\alpha}$ embedded curve. Suppose that $Y^2 \subseteq \H^{n+1}$ is a properly embedded minimal surface with a connected asymptotic boundary $\gamma$, then 
	\[
	\RA(Y) \leq - 2 \pi \chi(Y)
	\]	
	with equality if and only if $Y$ is a geodesic copy of $\H^2 \subseteq \H^{n+1}$ with boundary $\gamma \cong S^1$.
\end{corollary}
\noindent There is a similar formula for renormalized volume for $Y^4 \subseteq M^5$:
\begin{theorem}
Let $(M^5, g)$ be a Poincar\`e--Einstein manifold and $Y^4$ be a minimal hypersurface with $\partial Y = \partial M \cap \overline{Y}$. Then
\begin{align*}
6\RV(Y) &= 8 \pi^2 \chi(Y) - \int_Y \frac{|W|_Y^2}{4} dA_Y + \int_Y \frac{|E|_Y^2}{2} dA_Y \\
& \qquad - \int_Y \frac{|\hat{B}|^4}{24} dA_Y - \int \left[\frac{\Delta_Y(|\hat{B}|^2)}{2} + |\hat{B}|^2_Y \right] dA_Y
\end{align*}
where $W$ is the Weyl tensor, $E$ is the trace-free Ricci tensor, and $\hat{B}$ is the trace-free second fundamental form on $Y$.
\end{theorem}
\noindent which again provides connections to topological and conformally invariant quantities \cite{tyrrell2023renormalized}. At the moment, no rigidity theorem for higher dimensional renormalized volumes is known. In \cite{marx2021variations}, the author established regularity of minimal submanifolds near the boundary and computed variations of renormalized volume. This is described further in section \S \ref{MinSurfBackground}. \nl
\indent Physically, the AdS/CFT correspondence provides a major motivation to study inverse problems connected to renormalized area. Ryu and Takayanagi \cite{ryu2006aspects} proposed that for $A \subseteq \partial M$, the entanglement entropy is given by the renormalized area of the area-minimizing surface, $Y$, such that $\partial Y = \partial A$. The full nature of the AdS/CFT correspondence is unclear, as there are many pairs of boundary manifolds with a conformal class, $(N, [h])$, for which the bulk AdS manifold, $(M, g)$, is not unique. That being said, a first step towards the correspondence would be to show that the conformal field theory on the boundary manifold, $\partial M$, determines at least the asymptotic expansion near the boundary. In that the entanglement entropy \textit{should} be determined by the CFT, theorems \ref{ConformalInfinityThm} \ref{ExpansionThm} provides a partial affirmative answer to the question of determination: the entanglement entropy determines the asymptotic expansion of the bulk metric near $\partial M$. In some specific cases, e.g. when $g$ is log-analytic or $n+1 = 3$, we can establish uniqueness of the metric. 
%
%
\section{Statement of Results}
Let $(M^{n+1}, g)$ be an conformally compact metric. We state our results as follows
\begin{theorem} \label{ConformalInfinityThm}
Suppose that $g$ is partially even of order at least $2$, then the renormalized area on all minimal surfaces determines the conformal infinity, $c(g)$, of the ambient manifold.
\end{theorem}
\noindent The partially even condition is defined in \ref{partiallyEven} and is a slight generalization of an asymptotically hyperbolic metric with $\omega_1 = 0$. We remark that requiring $\omega_1 = 0$ is completely natural as renormalized volume loses its conformal invariance (and physical meaning) if $\omega_1 \neq 0$. See section \S \ref{PartiallyEvenAH} and the work of Bahaud--Mazzeo--Woolgar \cite{bahuaud2019ricci} for more motivation.
\begin{theorem} \label{ExpansionThm}
Let $(M^n, g)$, $g$ a partially even metric to order $2m^*$, and $n > 2m^*$. Suppose $\omega_0$ is known in equation \eqref{AHExpansion}. Then the renormalized volume on all $2m^*$-dimensional minimal submanifolds, $Y^m$, determines the asymptotic expansion in equation \eqref{AHExpansion} to arbitrary order.
\end{theorem}
\noindent \rmk \; Theorems \ref{ConformalInfinityThm} and \ref{ExpansionThm} fit together nicely if we restrict to $2$-dimensional minimal surfaces: one can first determine $c(g)$, and then after choosing a representative $\omega_0 \in c(g)$, the rest of the asymptotic expansion is determined. \nl \nl
In many ways, the ideas of theorems \ref{ConformalInfinityThm} and \ref{ExpansionThm} are inspired by the work of Graham--Guillarmou--Stefanov--Uhlmann \cite{graham2019x}. The authors consider a family of ``short geodesics" in section \S 4 which allows them to recover the metric from the renormalized length. We will see in section \S \ref{ProofSection} that recreating this technique for minimal surfaces which are close to geodesic hemispheres (i.e. ``short hemispheres") gives the \textit{trace} of the coefficients $\omega_i$ at a point $p$ on the boundary. In codimension $> 1$, this is enough to recover the asymptotic expansion of the metric. However, in the codimension $1$ case (which is of most interest), knowing the trace is not enough. To resolve this, we use variational techniques similar to those initiated in \cite{marx2021variations} to make full use of the ``short hemispheres" and recover the asymptotic expansion. \nl
\indent There are a few novel techniques compared to the geodesic case. Note that the renormalized length is \textit{not} conformally invariant, and depends on the choice of boundary defining function, or rather the associated representative of the conformal infinity - see equation 4.2 in \cite{graham2019x}. At the level of the boundary metric, the conformal \textit{invariance} of renormalized area means that the data of renormalized area on minimal surfaces is much coarser than the data of renormalized length on geodesics. To resolve this difficulty, we rely on the rigidity from corollary \ref{Rigidity} in order to detect an orthonormal frame for $\omega_0$. This allows us to determine the conformal infinity. We also note that we avoid Graham--Guillarmou--Stefanov--Uhlmann's use of the spherical tangent bundle as a means of parameterizing geodesics to second order. This clever idea is essential to their work, and while one could posit a similar parameterization of minimal hemispheres near a point on the boundary, we did not need to include this. Finally, many of the propositions and proofs in this paper rely on regularity results for minimal surfaces, as well as variational techniques for renormalized volume which have been developed in \cite{marx2021variations}. The author is unaware of an analogue of this in the setting of geodesics and renormalized length. \nl
\indent We can recreate the main theorems of Graham--Guillarmou--Stefanov--Uhlmann \cite{graham2019x} using renormalized volume. As immediate applications of theorems \ref{ConformalInfinityThm} and \ref{ExpansionThm}, we have the analogous boundary rigidity theorem for renormalized area:
\begin{corollary} \label{SurfaceGGSURigidity}
Suppose $(M, g)$, $(M, g')$ are two CC partially even metrics of order at least $2$. For $\gamma \subseteq \partial M$, let $Y^2_{\gamma, g}$ and $Y^2_{\gamma, g'}$ be the corresponding minimal surfaces (with respect to $g$ and $g'$ respectively) with $\partial Y^2_{\gamma, g} = \partial Y_{\gamma, g'}^2 = \gamma$. Suppose that $\RV_{g}(Y_{\gamma, g}^2) = \RV_{g'}(Y_{\gamma, g'}^2)$ for all $\gamma$. Then there exists $\psi: \overline{M} \to \overline{M}$ such that
\begin{align*}
	\psi \Big|_{\partial M} &= Id \\
	\psi^*(g') - g &= O(x^{\infty}) 
\end{align*}
Furthermore, if $g$ and $g'$ have log-analytic expansions and $\pi_1(\overline{M}, \partial \overline{M}) = 0$ then $\psi^*(g') = g$.
\end{corollary}
\noindent In section \S \ref{polyHomAnalytic}, we define ``log-analytic" formally but essentially it means that $g$ has polyhomogeneous expansion which is absolutely convergent in a collar neighborhood of the boundary. As in theorem \ref{ExpansionThm}, the same conclusion holds if we assume knowledge of the conformal infinity and a higher dimensional renormalized volume:
\begin{corollary} \label{GeneralGGSURigidity}
Suppose $(M, g)$, $(M, g')$ are two CC metrics that are partially even to order at least $m = 2m^*$ and $c(g) = c(g')$. For $\gamma^{m-1} \subseteq \partial M$, let $Y^m_{\gamma, g}$ and $Y^m_{\gamma, g'}$ be the corresponding minimal surfaces with $\partial Y^m_{\gamma, g} = \partial Y^m_{\gamma, g'} = \gamma$. Suppose that $\RV_{g}(Y_{\gamma, g}^m) = \RV_{g'}(Y_{\gamma, g'}^m)$ agree for all $\gamma$. Then there exists $\psi: \overline{M} \to \overline{M}$ such that
\begin{align*}
	\psi \Big|_{\partial M} &= Id \\
	\psi^*(g') - g &= O(x^{\infty}) 
\end{align*}
Furthermore, if $g$ and $g'$ have log-analytic expansions and $\pi_1(\overline{M}, \partial \overline{M}) = 0$ then $\psi^*(g') = g$.
\end{corollary}
%
%
\noindent We also apply our inverse problem to Poincar\`e-Einstein manifolds, $(M^{n+1}, g)$, for which 
\begin{align*}
(\text{n even}) \implies g &= \frac{dx^2 + \omega_0 + x^2 \omega_2 + (\text{even terms}) + x^{n} \omega_{n} + x^n \log(x) \omega_{n}^* + O(x^{n+1})}{x^2} \\
(\text{n odd}) \implies g &= \frac{dx^2 + \omega_0 + x^2 \omega_2 + (\text{even terms}) + x^{n-1} \omega_{n-1} + x^{n} \omega_{n} + O(x^{n+1})}{x^2}
\end{align*}
Moreover, $\omega_n$ ($n$ even) and $\omega_{n-1}$ ($n$ odd) are non-local terms and often thought of as the image of a Dirichlet-to-Neumann type map for Poincar\`e-Einstein metrics (see \S \ref{PartiallyEvenAH} for more background). With this, we state the following application
\begin{corollary} \label{DtoNDetermination}
Suppose $(M^{n+1},g)$ is a Poincar\`e-Einstein manifold and $\omega_0$ is known in equation \eqref{AHExpansion}. Then knowledge of the renormalized area on all minimal surfaces for any dimension $2m^* < n + 1$ determines the non-local coefficients, $\omega_{n}$ when $n$ even and $\omega_{n-1}$ when $n$ odd, in equations \eqref{PEExpansionEven} \eqref{PEExpansionOdd}.
\end{corollary}
\noindent In the particular case of $n = 2$, Poincar\`e-Einstein manifolds are known to be quotients of hyperbolic $3$-space by a convex, co-compact subgroup $\Gamma$. Let $\Sigma^2 = \partial M^3$, then we have the following global determination of the metric
\begin{corollary} \label{3DApp}
Suppose $(M^3,g)$ is a Poincar\`e-Einstein manifold with $\Sigma \cong \partial M$ and $\omega_0$ is known in equation \eqref{3DPEExpansion}. Then knowledge of the renormalized area on all minimal surfaces determines the conformal structure on $\Sigma$ and hence, $g$ globally.
\end{corollary}
\subsection{Acknowledgements}
The author wishes to thank
\begin{itemize}
\item Rafe Mazzeo for his encouragment, many helpful conversations, and a close reading of this paper.
\item Otis Chodosh for his advice on and suggestions for this project. 
\item Jeffrey Case and Steve Kerckhoff for providing key references for corollary \ref{3DApp}.
\item Aaron Tyrrell for many insightful conversations about renormalized volume and rigidity. 
\end{itemize}
The author is supported by a Stanford Graduate Fellowship and NSF Graduate Fellowship. The author dedicates this paper to his grandmother, Shirley Kuo.
\section{Background}
We present some background for the reader unfamiliar with renormalized area on asymptotically hyperbolic (AH) spaces. We refer the reader to \cite{marx2021variations}, \S 2 and \S 6 for a review of renormalized volume, boundary defining functions, and finite part evaluation through Riesz Regularization. See also \cite{albin2009renormalizing} \S 2.
\subsection{Conformally Compact, Asymptotically Hyperbolic manifolds}
\label{CCAH}
Let $(M^{n+1}, g)$ be a complete riemannian manifold. $M$ is \textit{conformally compact} if there exists a topological compactification, $\overline{M}$, along with $x: \overline{M} \to \R$, such that 
\begin{enumerate}
\item $\bg = x^2 g$ is a valid metric on all of $\overline{M}$
\item $x \Big|_{\partial M} = 0, ||\n^{\bg} x||_{\bg}^2 \Big|_{\partial M} > 0$
\end{enumerate} 
Such an $x$ is called a \textbf{boundary defining function}. Note that if $\varphi: \overline{M} \to \R^+$ is a positive smooth function, then $\tilde{\rho} = \varphi \rho$ is also a boundary defining function. As a result, we can consider the \textit{unique} conformal class $[\overline{g} |_{\partial M}]$ determined by $g$, called the \textit{conformal infinity}. With this choice of $x$, the metric $g$ will split as in equation \eqref{AHExpansion}. Such a splitting of the metric is motivated by Graham-Lee Normal Form \cite{graham1991einstein} \cite{fefferman2011ambient} for Poincar\'e-Einstein (PE) manifolds, which form a subclass of metrics that we mention for historical significance. These metrics have since been generalized to partially even AH metrics in \cite{bahuaud2019ricci}. \nl \nl
For PE manifolds with a chosen representative, $k_0$, in the conformal infinity, there exists a bdf $x$, for which $\overline{g}$ splits as in equation \eqref{AHExpansion} with further constraints on $\omega_k$. Moreover, $\omega(y,x)$ is regular up to order $x^m$ as shown in \cite{graham1991einstein}. The bdf $x$ is special if 
\[
||d \log(x)||_{g} = 1
\]
holds in a neighborhood of $\partial M$. Furthermore, by equation \eqref{AHExpansion} we have $\omega_0 = x^2 g \Big|_{\partial M}$. Given these conditions, $x$ is unique (see \cite{fefferman1985Conformal} for details). When $\omega_1 = \omega_3 = \cdots = \omega_{2m^*-1} = 0$ in equation \eqref{AHExpansion}, renormalized volume is conformally invariant for $2m^*$-dimensional minimal submanifolds of $M$ in the sense that it does not depend on the choice of $k_0 \in [\bg|_{\partial M}]$ and the corresponding special bdf used. Thus, we can define renormalized volume for even dimensional minimal submanifolds as long as we use a special bdf (see \cite{albin2009renormalizing} \cite{graham1999conformal}). 
\begin{example}
	Consider the Poincare Ball model of hyperbolic space $M = \H^3$. The metric on $\H^3$ is 
	\[
	g = \frac{4}{(1 - r^2)^2} \left[ dr^2 + r^2 d\phi^2 + r^2 \sin^2\phi d \theta^2 \right]
	\]
	is Einstein. We want to find a special bdf, $\rho$, for $\H^3$. We assume that it is rotationally symmetric, i.e. $\rho_{\theta} = \rho_{\phi} = 0$. With this, we compute 
	\[
	1 = ||d \log(\rho)||_{g}^2 = \frac{\rho_r^2}{\rho^2} g^{rr} = \partial_r(\log(\rho))^2 \frac{(1 - r^2)^2}{4}
	\]
	we take the negative root and get 
	\[
	\partial_r (\log(\rho)) = \frac{-2}{1 - r^2}
	\]
	Integrating and exponentiating, we compute
	\[
	\rho = A \frac{1 - r}{1 + r}
	\]
	for $0 \leq r \leq 1$ and some constant $A$. Note that as long as $A \neq 0$, we have $\rho^{-1}(0) = \{ r = 1\} = S^2$, which is the boundary of $\overline{\H^3}$. Suppose that we want to prescribe the standard metric on this boundary. i.e. $k_0(\theta) = \sin^2 \phi d \theta^2 + d \phi^2$. Then we have that 
	\[
	k_0 = \rho^2 h\Big|_{r = 1} = \frac{4 A^2}{(1 + r)^4} [dr^2 + r^2 d \phi^2 + r^2 \sin^2 \phi d \theta^2] \Big|_{r = 1} = \frac{4 A^2}{16} [d \phi^2 + \sin^2 \phi d \theta^2]
	\]
	so we choose $A = 2$ so that $\rho$ is positive. We see that $\bg \Big|_{r = 1} = \bg \Big|_{\partial M} = k_0$. Note that
	\begin{align*}
		\overline{\nabla} \rho & = \overline{g}^{ij} (\partial_i \rho) \partial_j = \overline{g}^{rr} (\partial_r \rho) \partial_r = \frac{(1 + r)^4}{16} \cdot \frac{-4}{(1 + r)^2} \partial_r\\
		\overline{\nabla} \rho \Big|_{r = 1} &=  - \partial_r
	\end{align*}
	which is non-zero.
\end{example}
\subsection{Partially Even Metrics} \label{PartiallyEvenAH}
In this section, we recall the notion of partially even metrics of order $m$, as set forth by Bahaud--Mazzeo--Woolgar in \cite{bahuaud2019ricci}. These metrics are natural to consider in that the renormalized volume of $m$-dimensional minimal submanifolds are necessarily conformally invariant. 
\begin{definition} \label{partiallyEven}
An AH metric $g$ is called even of order $2 \ell$, if in Graham--Lee normal form one has
\begin{equation} \label{EvenExpansion}
g = \frac{dx^2 + \omega_0 + x^2 \omega_2 + \dots + x^{2 \ell-2} \omega_{2 \ell-2} + O(x^{2 \ell} \log(x)^p)}{x^2}
\end{equation}
where $p > 0$ and $x$ a special boundary defining function.
\end{definition}
\noindent Note that this is a slight modification of the original definition 2.1 from \cite{bahuaud2019ricci}, as we will want to allow for $x^k \log(x)^p$ terms occuring at order $k = 2 \ell$ or higher. Often we will be concerned with the case of $(M^{n+1}, g)$ with $n$ even and there being a non-zero $x^n \log(x)$ term in our expansion (see \S \ref{Applications}). \nl \nl
\noindent An important result is that this parity condition is independent of special bdf, and also holds for bdfs as opposed to special bdfs
\begin{proposition}[Bahaud--Mazzeo--Woolgar, Proposition 2.5]
If $g$ is a partially even metric of order $2 \ell$ with respect to $x$, a special bdf, then it is partially even with respect to any other special bdf, $x'$.
\end{proposition}
\begin{lemma}[Bahaud--Mazzeo--Woolgar, Lemma 2.8]
$g$ is a partially even metric of order $2 \ell$ with respect to any boundary defining function if and only if it is a partially even metric of order $2 \ell$ with respect to any special boundary defining function
\end{lemma}
\noindent \Pf The proof of both statements (with the $\log(x)$ modifications) follow the exact same procedure as in the original proofs from \cite{bahuaud2019ricci} \S 2.2, taking slight care to handle the $x^{2 \ell} \log(x)^p$ terms. \qed \nl \nl
%
We now describe the main example of partially even metrics, Poincar\'e--Einstein metrics. Recall that Poincar\'e-Einstein metrics are slight modifications of conformal compact asymptotically hyperbolic metrics. For $(M^{n+1}, g)$ with $n$ even,
\begin{align} \label{EinsteinEq}
	\Ric &= -n g \\ \label{PEExpansionEven}
	g &= \frac{dx^2 + \omega_0 + x^2 \omega_2 + (\text{even terms}) + x^{n} \omega_{n} + x^n \log(x) \omega_{n}^* + O(x^{n+1})}{x^2}
\end{align}
In general, the expansion for $g$ is even up to order $n$, and then loses parity and allows for $x^k \log(x)$ terms (i.e. a \textit{polyhomogeneous} expansion).  Moreover, $\omega_{2k}$ and $\omega_n^*$ are determined by $\omega_0$ for all $k \leq (n-2)/2$. $\omega_n$ is traceless but otherwise undetermined and knowledge of $\omega_0, \omega_n$ determines the expansion of $g$ to all orders. When $n$ is odd, we have
\begin{align} \label{PEExpansionOdd}
	g &= \frac{dx^2 + \omega_0 + x^2 \omega_2 + (\text{even terms}) + x^{n-1} \omega_{n-1} + x^{n} \omega_{n} + O(x^{n+1})}{x^2}
\end{align}
in particular, there are no logarithm terms and the expansion is AH. In this case, $\omega_{n-1}$ is the undetermined term and all lower order terms are determined by $\omega_0$. Equations \eqref{PEExpansionEven} \eqref{PEExpansionOdd} are consequences of equation \eqref{EinsteinEq} as shown by Graham-Lee in their seminal work \cite{graham1991einstein}, and is titled the Graham--Lee normal form. See also \cite{albin2009renormalizing} section $2$ for more background.
\subsection{Log-analytic expansions} \label{polyHomAnalytic}
Let $x$ be a boundary defining function for $\overline{M}$, and $s$ a coordinate function on the boundary. We define the space of conormal functions on $\overline{M}$, $\mathcal{A}(\overline{M})$, to be 
\[
f(s,x) \in \mathcal{A}(\overline{M}) \implies \forall j, \beta, \; \exists \; C_{j,\beta} \; \st \; |(x \p_x)^j \p_s^{\beta} f| \leq C_{j, \beta}
\]
In this paper we will be concerned with expansions of the form
\begin{equation} \label{polyHomExpansion}
T(s,x) = \sum_{k = 0}^{M} x^k \sum_{i = 0}^{N_k} T_{ki}(s) \log(x)^i + O(x^{M+1} \log(x)^{N_{M+1}})
\end{equation}
functions obeying the expansion in equation \eqref{polyHomExpansion} represent a subset of \textit{polyhomogeneous} functions. We will say that a function with the above expansion is \textbf{polyhomogeneous} if
\[
T(s,x) - \sum_{k = 0}^{M} x^k \sum_{i = 0}^{N_k} T_{ki}(s) \log(x)^i \in x^{M+1} \mathcal{A}
\]
We will then say that $T(s,x)$ is \textbf{log-analytic} on a neighborhood of the boundary if $x$ is a real-analytic defining function for $\partial M$ and $T(s,x)$ admits an expansion as in equation \eqref{polyHomExpansion} which converges absolutely in $x$ for $|x| < r$ for some $r > 0$.
%
\subsection{Renormalized Volume through Riesz Regularization}
\label{Riesz}
Having defined special bdfs, we can define renormalized volume via Riesz regularization: given a partially even metric of order $m$ even, $g$, let $x$ be a special bdf
\begin{align*}
	f&: \{ \text{Re}(z) > m \} \to \C \\
	f(z) &= \int_Y x^z dA_Y 
\end{align*}
$f(z)$ is holomorphic for $\text{Re}(z) > m$, and can be extended to be meromorphic on $\C$ with poles at $ z \in \{- \infty, \dots, -1, 0, 1, \dots, m\}$ (see \cite{melrose1996homology}\cite{paycha2003heat} for a proof). We define
\[
\RV(Y^m):= \FPz \int_Y x^{z} dA_Y = \FPz f(z)
\]
Computing $\FPz f(z)$ amounts to subtracting off the pole at $z = 0$ (if it exists) and evaluating the remaining difference. This process is known as Riesz regularization and see \cite{marx2021variations}, \S 7 for more details.
\begin{example}
	We compute $\RV(\H^2)$ for $\H^2 \subseteq \H^3$ using Riesz regularization in the half space model (we leave it to the reader to compute this for the poincare ball model).
	\[
	\zeta(z) = \int_{\H^2} x^z dA_{\H^2} = \int_{x = 0}^{1}\int_{\theta = 0}^{2\pi} x^{z - 2} dx d\theta = 2\pi \left[\frac{x^{z-1}}{z - 1}\right]_{x = 0}^{x = 1} = 2 \pi \frac{1}{z - 1}
	\]
	Again, when we find the meromorphic extension, we first assume $\text{Re}(z) \gg 0$ so that $0^{z-1} = 0$. There is no pole at $z = 0$ in this extension, so 
	\[
	\RV(\H^2) = \FPz \zeta(z) = \zeta(0) = \boxed{-2 \pi}
	\]
\end{example}

\subsection{Minimal Surfaces in AH spaces} \label{MinSurfBackground}
In \cite{marx2021variations}, the author expanded upon the study of minimal surfaces in hyperbolic space, extending pre-existing regularity results. Much of this work was based in the work of Alexakis--Mazzeo \cite{alexakis2010renormalized} and we recall the following background. \nl \nl
Given $(M, g)$ conformally compact asymptotically hyperbolic, consider $x: \overline{M} \rightarrow \R^{\geq 0}$ a special boundary defining function, along with  $\gamma^1 \subseteq \partial M$, a curve lying in the boundary. We define the \textbf{boundary cylinder} submanifold
\[
\Gamma = \{(s,x) \; | \; s \in \gamma, \;\; x \geq 0\}
\]
We can now consider fermi coordinates off of this submanifold via 
\[
(x,s,z) \sim \overline{\exp}_{\Gamma, (s,x)}(z \overline{N})
\]
where $\overline{\exp}$ and $\overline{N}$ are the exponential map and unit normal vector to $\Gamma$ with respect to the compactified metric, $\bg$. We then have the following from \cite{marx2021variations}
\begin{theorem*}[MK, Theorem 3.1]
Let $\gamma$ be a $C^{m,\alpha}$ embedded curve in $\partial M$ where $(M^{n+1}, g)$ is Poincar\`e-Einstein, $k \geq m$, $0 < \alpha < 1$. Let $Y^m$ be an even minimal surface with asymptotic boundary $\partial Y = \gamma$. Then $Y$ can be described graphically over $\Gamma$ as 
\[
Y = \{ (s,z = u(s,x), x) \}
\]
for $u$ even to order $m$ in $x$, i.e.
\begin{equation} \label{UExpansion}
u(s,x) = u_2(s) x^2 + u_4(s) x^4 + \dots + u_m(s) x^m + u_{m+1}(s) x^{m+1} + O(x^{m+2})
\end{equation}
\end{theorem*}
\noindent This is a more general version of the expansion given for $Y^2 \subseteq M^3$ in \cite{alexakis2010renormalized}.
In particular, $u_{m+1}(s)$ is non-local, i.e. cannot be determined by $\gamma$, and any higher order terms in equation \eqref{UExpansion} are determined by $u_{m+1}(s)$ and $\gamma$. In \cite{marx2021variations}, the author describes how $\gamma \mapsto u_{m+1}(s)$ could be thought of as a Dirichlet-to-Neumann map for minimal surfaces in asymptotically hyperbolic spaces, analogous to the map $g_{(0)} \rightarrow g_{(n)}$ considered in the ambient PE case (see equation 1.9, \cite{anderson2005geometric}). \nl \nl
Given $Y^m \subseteq M^{n+1}$, we consider $h = g\Big|_{TY}$ the induced metric and $\bh = \bg \Big|_{TY}$ the induced compactified metric. If $s$ is a local coordinate for $\gamma$ and $x$ our special bdf, then near the boundary we can write
\begin{equation} \label{AreaFormEqn}
dA_{Y_{\gamma, \delta}}(s,x) = \frac{1}{x^m}\sqrt{\det \bh} dx dA_{\gamma(y)}(s) 
\end{equation}
where $\sqrt{\det \bh}$ is even up to order $m$. We also recall the following result about the normal vector field to $\nu_Y$: let $\{s_a\}_{a = 1}^{m-1}$ be coordinates for $\gamma^{m-1}$ and $z$ a normal coordinate to $\gamma^{m-1} \subseteq (\partial M)^m$. 
\begin{align} \label{NormalExpansion}
\nu &= x \bnu \\ \nonumber
\bnu &= c^z \partial_z + c^x \partial_x + c^a \partial_{s_a} \\ \nonumber
c^z &= 1 + \tilde{R}_z,  \\ \nonumber
c^x &= -u_x + \tilde{R}_x  \\ \nonumber
c^a &= -u_a + \tilde{R}_z 
\end{align}
where $\tilde{R}_z = O(x^2)$ and is even to order $m+2$, $\tilde{R}_x = O(x^3)$ and is odd to order $m+2$, and $\tilde{R}_z$ is even to order $m+2$. See \cite{marx2021variations} \S 10 for more details. Finally, if $T$ is a tensor with the analogous parity constraints, i.e. suppose
\[
T = T_{xx} dx^2 + T_{ax} ds^a dx + T_{ab} ds^a ds^b
\]
And $T_{xx}$ is even to order $m$, $T_{ax}$ is odd to order $m-1$, and $T_{ab}$ is even to order $m$. Then clearly
\[
T(\nu, \cdot) = A^x dx + B^a ds^a
\]
with $A^x$ being odd to order $m-1$ and $B^a$ even to order $m$.
\subsection{Background on Inverse Problems}
Here we describe past research in the area of inverse problems which were the main inspiration for this work. \nl
\indent In \cite{graham2019x}, Graham-Guillarmou-Stefanov-Uhlmann consider a boundary rigidity problem for asymptotically hyperbolic (AH) metrics via the $X$-ray transform. As part of this, they consider the renormalized length functional
\[
L(\gamma) = \lim_{\eps \to 0} \left( \ell_g(\gamma \cap \{\rho > \eps\}) + 2 \log(\eps) \right)
\]
and demonstrate that knowledge of the renormalized length functional on a family of ``short geodesics" determines the asymptotic expansion to any order.
\begin{theorem}
Let $\overline{M}$ be a compact connected manifold-with-boundary and let $g$, $g'$ be two asymptotically hyperbolic metrics on $M$. Suppose for some choices $h$ and $h'$ of conformal representatives in the conformal infinities of $g$ and $g'$, the renormalized lengths agree for the two metrices, i.e. $L_g = L_g'$. Then there exists a diffeomorphism $\psi: \overline{M} \to \overline{M}$ which is the identity on $\partial \overline{M}$ such that $\psi^* g' - g = O(\rho^{\infty})$ at $\partial \overline{M}$.
\end{theorem}
\noindent The idea of ``short geodesics" are essential in their theorem and we use a similar idea of ``short hemispheres" in the proof of our theorem \ref{ExpansionThm}. \nl 
\indent Another influential work is that of Alexakis--Balehowsky--Nachman \cite{alexakis2020determining} who describe a similar inverse problem but in the Riemannian setting. The authors consider $(M^3,g)$ a $C^4$ Riemannian manifold homeomorphic to $B^3$, along with a foliation of minimal surfaces from all directions and some geometric assumptions. In this case, knowledge of the area of minimal surfaces determines the metric.
\begin{theorem}[Alexakis-Balehowksy-Nachman, Theorem 1.4]
Let $(M, g)$ be homeomorphic to a $3$-ball with mean convex boundary which is $\eps_0$-$C^3$ close to the euclidean metric. Further assume that there is a foliation of $\partial M$ by simple closed curves, $\{\gamma(t)\}$, which induces a foliation $\{Y(t)\} = M$ by properly embedded, area minimizing surfaces with some regularity. Suppose that for each $\gamma(t)$ and any nearby perturbation, $\gamma(s,t) \subset \partial M$, we know the area of the properly embedded area minimizing surface $Y(s,t)$ with $\partial Y(s,t) = \gamma(s,t)$. Then the metric $g$ is uniquely determined up to isometries which fix the boundary.
\end{theorem}
\noindent Their proof relies on rigidity results for the Dirichlet-to-Neumann map for the Jacobi operator on any $Y(t)$. This differs from the techniques used in the proof of theorem \ref{ExpansionThm}, however it would be interesting to pursue an asymptotically hyperbolic version of metric rigidity through knowledge of the non-local term, $u_{m+1}(s)$, for all minimal surfaces. We also mention the recent work of Carstea--Lassas--Liimatainen--Tzou \cite{carstea2023inverse} who demonstrate rigidity of the surface itself, based on a graphical Dirichlet-to-Neumann map.
\section{Preliminaries and Notation}
\subsection{Dilation Maps}
In this section, we define a dilation map which will be used throughout the paper. Given $p \in \partial M$, we let $y$ be a coordinate system in a neighborhood $U(p) \subseteq \partial M$ such that $y = 0 \leftrightarrow p$. We then have on $(x,y) \in [0, \kappa) \times U(p)$ for some $\kappa$ sufficiently small, we can expand our partially even metric of order $2m^*$.
\begin{align*}
g & = g(x,y) \\
&= \frac{dx^2 + \sum_{k = 0}^{m^*-1} x^{2k} \omega_{2k} + \sum_{k = 2m^*}^M \sum_{\ell = 0}^P x^{k} \log(x)^{\ell} \omega_{k \ell} + O(x^{M+1} \log(x)^{P+1})}{x^2} \\
&= \frac{dx^2 + \sum_{k = 0}^{m^*-1} x^{2k} \omega^{2k}_{ij}(y) dy^i dy^j + \sum_{k = 2m^*}^M \sum_{\ell = 0}^P x^{k} \log(x)^{\ell} \omega^{k \ell}_{ij} dy^i dy^j + O(x^{M+1} \log(x)^{P+1})}{x^2}
\end{align*}
i.e. the metric has a polyhomogeneous expansion with log terms only appearing at orders $k \geq 2m^*$. As indicated above, we will adopt the following index notation:
\begin{align*}
\omega_k &:= \omega^k_{ij} dy^i dy^j \\
\omega_{k \ell} &:= \omega^{k \ell}_{ij} dy^i dy^j 
\end{align*}
Let 
\begin{align} \label{DeltaDiffeo}
F_{\delta} &: [0, \delta^{-1} \kappa) \times B_{\delta^{-1}}^{n}(0) \rightarrow \overline{M} \\
F_{\delta}(\tx, \ty) &= (\delta \tx, \delta \ty)
\end{align}
having implicitly identified a point $(x,y) \in [0,\kappa) \times U(p)$ with its corresponding point in $\overline{M}$. Suppose that $g$ is partially even up to order $2m^*$. Then we compute
\begin{align} \label{gDeltaMetric}
g_{\delta} &:= F_{\delta}^*(g) \\ \nonumber
&= \frac{d\tx^2 + \sum_{k = 0}^{m^*-1} \delta^{2k} \tx^{2k} \omega_{2k}(\delta \ty)}{\tx^2} \\ \nonumber
& + \frac{\sum_{k = 2m^*}^M \sum_{\ell = 0}^{P_k} \delta^k \tx^{k} [\log(\delta) + \log(\tx)]^{\ell} \omega_{k\ell}(\delta \ty) + O(\delta^{M+1} \log(\delta)^{P_{M+1}+1})}{\tx^2} \\ \nonumber
&= \frac{d\tx^2 + \sum_{k = 0}^{m^*-1} \delta^{2k} \tx^{2k} \omega_{2k}(\delta \ty)}{\tx^2} \\ \nonumber
& + \frac{\sum_{k = 2m^*}^M \delta^k \tx^{k} \sum_{\ell = 0}^{P_k} \log(\delta)^{\ell} \sum_{j = 0}^{P_k - \ell} \log(\tx)^j \begin{pmatrix} \ell + j \\ \ell \end{pmatrix} \omega_{k (\ell + j)} + O(\delta^{M+1} \log(\delta)^{P_{M+1}+1})}{\tx^2} 
\end{align}
One can compute 
\begin{align} \label{NthMetricDerivNoLog}
2N \leq 2m^* - 2, \qquad \frac{d^{2N}}{d \delta^{2N}} g_{\delta} \Big|_{\delta = 0}(\tx, \ty) &= \sum_{k = 0}^{N} c_{2k} \tx^{2k} (D^{2N-k} \omega_k)(\ty^{2N-k}) \\ \nonumber
&= (D^{2N} \omega_0)(\ty^{2N}) + \dots + (2N!)\tx^{2N} \omega_{2N}(0)
\end{align}
where $c_{2k} > 0$ is some combinatorial constant. Recall that $y = 0 \leftrightarrow \ty = 0 \leftrightarrow p$, so the $N$th derivative of $g_{\delta}$ depends on $\omega_N(p)$. Of particular importance is that if we know $\omega_0, \dots, \omega_{2N-2}$ and $\p_{\delta}^{2N} g_{\delta}$, then we can recover $\omega_{2N}$. For larger values of $N$ (not necessarily even), we have to account for the potential $\log(\delta)$ terms. We first note that
\[
\frac{1}{\log(\delta)^a} \frac{d^b}{d \delta^b} \delta^{c} \log(\delta)^d \Big|_{\delta = 0} = \begin{cases}
	0 & c > b \\
	b! & a = d, \; c = b \\
	0 & c = b,  \; d < a
	\end{cases}
\]
Immediately from this, 
if we let $\FPd$ denote computing the finite part as $\delta \to 0$ (i.e. removing all factors of $\{\delta^{-k}\log(\delta)^p\}_{k, p > 0}$), then we also have
\begin{align} \label{FPDelta}
\FPd \frac{1}{\log(\delta)^a} \frac{d^b}{d \delta^b} \delta^{c} \log(\delta)^d \Big|_{\delta = 0} &= \begin{cases}
	b! & d = a, \; c = b \\
	0 & d < a, \; c = b \\
	0 & c \neq b 
\end{cases} \\ \nonumber
\implies \FPd \frac{1}{\log(\delta)^a} \frac{d^b}{d \delta^b} g_{\delta} &= b! \tx^{b-2} \omega_{ba} + F(\{\tx^{b-2} \omega_{b \ell}\}_{\ell > a})
\end{align}
%
where $F$ is a computable linear function in its entries. For all $\delta < \delta_0$ sufficiently small, we will often consider $g_{\delta}$ restricted to 
\begin{align} \label{OpenSets}
V_0 &= [0, 100) \times B_{100}(0) \subseteq \R^{\geq 0} \times \R^n \\ \nonumber
V_1 &= [0, 75) \times B_{75}(0) \\ \nonumber
V_2 &= [0, 50) \times B_{50}(0)
\end{align}
\subsection{Parameterization of hemispheres in hyperbolic space} \label{HemisphereParamSection}
Let $(x, y_1, \dots, y_n) \in \H^{n+1}$ with the standard half-space metric. We recall the parameterization of the hemisphere, $HS^m$
\begin{align*}
HS^m & = \{(x, y_1, \dots, y_{m}, y_{m+1}, \dots, y_n) = (x, \sqrt{1 - x^2} F_{m-1}(\vec{\theta}), y_{m+1} = 0, \dots, y_n = 0) \}
\end{align*}
here, $F_{m-1}(\vec{\theta}) \in S^{m-1}$ is the standard polar coordinate parameterization of the $m-1$ dimensional sphere. As a result, we have the following basis for the tangent space and induced metric coefficients
\begin{align} \label{HSMetric}
v_x &= \p_x - \frac{x}{\sqrt{1 - x^2}} F_m(\theta) \\ \nonumber
v_{\theta_i} &= \sqrt{1 - x^2} F_{m, i}(\theta) \\ \nonumber
h_{xx} &= \frac{1}{(1 - x^2)x^2} \\ \nonumber
h_{i j} &= g(v_{\theta_i}, v_{\theta_j}) = \frac{1 - x^2}{x^2} \bh_{ij} \\ \nonumber
h_{x i} &= g(v_x, v_{\theta_i}) = 0 \\ \nonumber
g(\nabla_{v_x} v_x, v_x) &= \frac{1}{2} v_x h_{xx}  \\ \nonumber
&= -\frac{1 - 2x^2}{(1 - x^2)^2 x^3} \\ \nonumber
g(\nabla_{v_x} v_x, v_{\theta_i}) &= - \frac{1}{2} v_{\theta_i} h_{xx} \\ \nonumber
&= 0 \\ \nonumber
g(\nabla_{v_{\theta_i}} v_{\theta_j}, v_{\theta_k}) &= \frac{1}{2} [v_{\theta_i} h_{jk} + v_{\theta_j} h_{ik} - v_{\theta_k} h_{ij}] \\ \nonumber
&= \frac{1 - x^2}{x^2} \bGamma_{ijk}\\ \nonumber
g(\n_{v_{\theta_i}} v_{\theta_j}, v_x) &= - \frac{1}{2} v_x h_{ij}  \\ \nonumber
&= \frac{1}{x^3} \bh_{ij}
\end{align}
where $\bh_{ij}$ and $\bGamma_{ijk}$ are the metric coefficients and Christoffel symbols on $S^m$ with respect to the euclidean metric. As a consequence, we have 
\[
dA_{HS^m} = \frac{(1-x^2)^{m-2}}{x^{2m}} dA_{S^{m-1}} dx
\]

\subsection{Mechanics of Finite Part Evaluation} \label{mechanics}
This is a summary of computing the finite part of integrals with respect to a complex parameter, $z$. See \cite{marx2021variations} \S 6 for more details. \nl 
\indent When computing variations of renormalized volume, we encounter integrals of the form 
\[
I(z) = \int_Y z^p b(x,s) x^{z-j} dA_Y
\]
for $b(x,s)$ having a polyhomogeneous expansion in $x$ (after pulling back to $\Gamma$) and $i \geq 0$, $j \in \{0,1,2\}$, and $p \in \{0,1\}$. We write 
\[
I(z) = \left(\int_{Y \cap \{x < \eta\}} + \int_{Y \cap \{x \geq \eta\}}\right) z^p b(x,s) x^{z-j} dA_Y = I_1(z, \eta) + z^p I_2(z, \eta)
\]
for some $1 \gg \eta > 0$, where we've pulled out the factor of $z^p$ in the $\{x \geq \eta\}$. As before, $I_2(z)$ is holomorphic because the integral is over $x \geq \eta$. In particular
\[
\FPz z^p I_2(z, \eta) = \begin{cases}
	I_2(0, \eta) & p = 0\\
	0 & p = 1
\end{cases}
\]
We further assume the following expansions (after pulling back to Fermi coordinates)
\begin{align}  \nonumber
	dA_Y &= \frac{\sqrt{\det \overline{h}}}{x^m} dx \wedge dA_{\gamma} \\ \label{AreaExpansion}
	\sqrt{\det \overline{h}} &= \sum_{k = 0}^{m+2} \overline{q}_j(s) x^k + \tilde{Q}(s) x^m \log(x) + \overline{Q}(s) x^{m+1} \log(x) + O(x^{m+2} \log(x)) \\ \label{bExpansion}
	b(x,s) &= \sum_{k = 0}^{m+2} b_j(s) x^k + \tilde{B}(s)x^m \log(x) + B(s) x^{m+1} \log(x) + O(x^{m+2} \log(x))
\end{align}
i.e. if a $x^d \log(x)^q$ term manifests, it can occur only when $d \geq m$. $I_1$ expands as 
\begin{align*}
	I_1(z, \eta) &= z^p \int_{0}^{\eta} \int_{\gamma} x^{z-m-j}\left[\sum_{\ell = 0}^{m+2} \sum_{k + j = \ell} b_{k} \overline{q}_j x^{\ell} \right] ds dx  \\
	& + z^p \int_{0}^{\eta} \int_{\gamma} x^{z-m-j}\left[ (\overline{q}_0 \tilde{B} + \tilde{Q} b_0 ) x^m \log(x) + (b_0 \overline{Q} + b_1 \tilde{Q} + \overline{q}_0 B + \overline{q}_1 \tilde{B}) x^{m+1} \log(x)\right] ds dx  \\
	&+ z^p \int_{0}^{\eta} \int_{\gamma} x^{z-m-j}\left[ (b_1 \overline{Q} + b_2 \tilde{Q} + \overline{q}_1 B + \tilde{B} \overline{q}_2)x^{m+2} \log(x) \right] ds dx \\
	& + z^p\int_0^{\eta} \int_{\gamma} O(x \log(x)) ds dx
\end{align*}
Observe that
\[
\FPz z^p \int_{0}^{\eta} \int_{\gamma} O(x \log(x)) dx = \begin{cases}
	C(\eta) & p = 0 \\
	0 & p = 1
\end{cases}
\]
for some finite constant $C(\eta)$. It remains to compute
\begin{align} \label{FPExpansion}
	\FPz I_1(z, \eta) &= \FPz z^p \sum_{k = 0}^{m+2} c_k \int_{0}^{\eta} x^{z+ k - m - j} dx  \\ \nonumber
	& + \FPz z^p \sum_{k = 0}^2 c_{m + k}^* \int_0^{\eta} x^{z + k - m - j} \log(x) dx
\end{align}
for 
\begin{align} \label{ckCoeffs}
	c_k &= \int_{\gamma} \left[\sum_{\ell + j = k} b_{\ell} \overline{q}_j  \right] dV_{\gamma}, \qquad 0 \leq k \leq m + 2 \\ \nonumber
	c_m^* &= \int_{\gamma} \left[\overline{q}_0 \tilde{B} + \tilde{Q} b_0 \right] dA_{\gamma} \\
	c_{m+1}^* &= \int_{\gamma} \left[ b_0 \overline{Q} + b_1 \tilde{Q} + \overline{q}_0 B + \overline{q}_1 \tilde{B}\right] dA_{\gamma} \\ \nonumber
	c_{m+2}^* &= \int_{\gamma} \left[ b_1 \overline{Q} + b_2 \tilde{Q} + \overline{q}_1 B + \tilde{B} \overline{q}_2\right] dA_{\gamma}
\end{align}
Integrating,
\[
I_1(z, \eta) = z^p \left(c_{m+j-1}\frac{\eta^z}{z} + c^*_{m+j-1}\frac{\eta^z ((z \log(\eta) - 1)}{z^2}+ F(\eta, z)\right) 
\]
where $F(\eta, z)$ is a holomorphic near $0$ and computable in terms of the remaining $c_k$ coefficients through equation \eqref{FPExpansion}. In particular
\begin{align*}
	\FPz z^p F(\eta, z) &= \begin{cases}
		F(\eta, 0) & p = 0 \\
		0 & p = 1
	\end{cases} \\
	\FPz z^p \frac{\eta^z}{z} &= 
	\begin{cases}
		\log(\eta) & p = 0\\
		1 & p = 1 
	\end{cases} \\
	\FPz z^p \frac{\eta^z (z \log(\eta) - 1)}{z^2} &=
	\begin{cases}
		\frac{\log(\eta)^2}{2} & p = 0 \\
		0 & p = 1 
	\end{cases}
\end{align*}
We summarize this as
\begin{lemma} \label{FiniteEvaluation}
	Consider integrals of the form
	\[
	I(z) = \int_Y z^p b(x,s) x^{z-j} dA_Y 
	\]
	for $b$ and $dA_Y$ having polyhomogeneous expansions in $x$ and $p \in \{0,1\}$, $j \in \{0,1,2\}$. Moreover, assume that $x^d \log(x)^q$ terms only appear when $d \geq m$ and $q = 1$, or $d \geq m + 3$. Then we have that 

\[
		\FPz I(z) = \begin{cases} 
			C(\eta) + F(\eta, 0) + I_2(0, \eta) + c_{m+j-1}\log(\eta) + c^*_{m+j-1}\frac{\log(\eta)^2}{2} & p = 0 \\[1ex]
			c_{m+j-1} & p = 1 
		\end{cases}
\]
	for the coefficients $\{c_{k}\}$ listed above.
\end{lemma}
\noindent \rmk:
\begin{itemize}
	\item This calculation illustrates the following key point: \textbf{when $p = 1$, the finite part can be expressed as an integral over the boundary}. We will refer to this process from here on as \textbf{localization}. 
	\item Taking $b(x,s) = 1$ and $p = 0$ demonstrates how to compute the renormalized volume of $Y$ via Riesz regularization
	
	\item While the result for $p = 0$ seems to depend on $\eta$, one can show that by changing $\eta \to \eta'$ and keeping track of boudary terms from the intermediate integral $\int_{x = \eta}^{\eta'}$, the result is independent of $\eta$. This is shown in \cite{marx2021variations}.
	
	\item Often, the coefficients $\{b_j(s), B(s), \tilde{B}(s), \overline{q}_j(s), \overline{Q}(s), \tilde{Q}(s)\}$ from equations \eqref{AreaExpansion} \eqref{bExpansion} will also be functions of an external parameter, $\delta$. When the $\delta$-dependence is $C^k$ or polyhomogeneous of the correct order, one can interchange differentiation and other operations involving $\delta$ (e.g. removing finite parts as $\delta \to 0$) with computing $\FPz$ in $x$. This is because computing $\FPz$ does not affect the boundary integrals used to compute the $c_k$ coefficients in equation \eqref{ckCoeffs}, which appear in the computation of $\FPz$ in lemma \ref{FiniteEvaluation}.
\end{itemize}

\section{Existence of short minimal surfaces}
Recall the open sets defined in equation \eqref{OpenSets}
\begin{align*} 
V_0 &= [0, 100) \times B_{100}(0) \subseteq \R^{\geq 0} \times \R^n \\ 
V_1 &= [0, 75) \times B_{75}(0) \\ 
V_2 &= [0, 50) \times B_{50}(0)
\end{align*}
We establish the following proposition about the existence of minimal surfaces with respect to $g_{\delta}$.
\begin{proposition} \label{ExistenceOfSmall}
Let $(M^{n+1}, g)$ be a partially even, asymptotically hyperbolic metric of order at least $2m+2$. Suppose $\gamma^{2m-1}$ is a $C^{2m,\alpha}$ embedded curve inside of $B_{50}(0) \subseteq \R^n$ with $Y^{2m}_{\gamma, 0}$ a non-degenerate minimal submanifold with respect to $g_{\delta = 0} \cong g_{\H^{n+1}}$ lying inside $V_1$. There exists a $\delta_0 = \delta_0(\gamma)$ such that for all $\delta < \delta_0(\gamma)$, there exists a submanifold, $Y^{2m}_{\gamma, \delta} \subseteq V_0$, minimal with respect to $g_{\delta}$ and asymptotic to $\gamma$.
\end{proposition}
\noindent \Pf We refer to proposition 4.1 in Alexakis--Mazzeo \cite{alexakis2010renormalized} for the existence of $Y_{\gamma, 0}$ and note that we use their same notion of non-degeneracy pertaining to the Jacobi operator on weighted spaces
\[
J_Y: x^{\mu} \Lambda_0^{2,\alpha} \rightarrow x^{\mu} \Lambda_0^{0,\alpha}
\]
for $0 < \alpha < 1$ and $-1 < \mu < m$. Here, the choice of weighting, $\mu$, guarantees invertibility. To show the existence of $Y_{\gamma, \delta}$, we search for graphical solutions expressed through the exponential map of the form
\begin{align*}
Y_{\gamma, \delta} &= \exp_{Y_{\gamma, 0}, g_0}(\dot{\phi}_{\delta}(p) \nu(p)) 
\end{align*}
And compute the linearization of the mean curvature equation
\begin{align} \nonumber
\frac{d}{d \delta} H(Y_{\gamma, \delta}, g_{\delta}) &= J_Y(\dot{\phi}) + L_2(\dot{g})  \\
\dot{g} &= \frac{\nabla \omega^0\Big|_P(\ty) \cdot \ty}{\tx^2} \\ \label{L2Eqn}
L_2(\dot{g}) & = - \langle \dot{g}, A_{Y_{\gamma}} \rangle + \div_{Y_{\gamma}}(\dot{g}(\nu, \cdot)) - \frac{1}{2} \nu(\tr_{Y_{\gamma}}(\dot{g}))
\end{align}
where $Y_{\gamma} = Y_{\gamma,0}$ and $\nu = \nu_{Y_{\gamma}}$ and $\dot{g}$ denotes the derivative of $g_{\delta}$ (as in equation \eqref{gDeltaMetric}) at $\delta = 0$. Define
\begin{align*}
F_{\delta}&: x^{\mu} \Lambda^{2,\alpha}_0 \to x^{\mu} \Lambda^{2,\alpha}_0 \\
F_{\delta}(u) & = -\delta^{-1} J_{Y_{\gamma, 0}}^{-1}(H(Y_{\gamma}(\delta u), g_{\delta})) + u
\end{align*}
Here, $Y_{\gamma}(\delta u)$ denotes the image of $\exp_{Y_\gamma, g_0}(\delta u \nu)$ where $u: Y \to \R$, and $H(Y, g_{\delta})$ denotes the scalar mean curvature of $Y$ with respect to $g_{\delta} \Big|_{TY}$. We note that 
\begin{align*}
F_{\delta}(u) &= u \\
\implies H(Y_{\gamma}(\delta u), g_{\delta}) &= 0
\end{align*}
We compute
\begin{align*}
F_{\delta}(u) - F_{\delta}(v) &= -\delta^{-1}J_{Y_{\gamma, 0}}^{-1}[H(Y_{\gamma}(\delta u), g_{\delta}) - H(Y_{\gamma}(\delta v), g_{\delta})] + (u - v)
\end{align*}
Note that
\begin{align} \label{MCExpansionEqn}
H(Y_{\gamma}(\delta u), g_{\delta}) &= H(Y_{\gamma}(0), g_0) + [J_{Y_{\gamma,0}}(u) + L_2(\dot{g})] \delta + \overline{L}_{\delta}(u) + Q_{\delta}(u) + W_{\delta} \\ \nonumber
&= [J_{Y_{\gamma,0}}(u) + L_2(\dot{g})] \delta + \overline{L}_{\delta}(u) + Q_{\delta}(u) + W_{\delta} \\ \nonumber
||\overline{L}_{\delta}(u)|| & \leq K \delta^2 ||u|| \\ \nonumber
||Q_{\delta}(\delta u) || & \leq K ||\delta u||^2 \\ \nonumber
||W_{\delta}|| & \leq K \delta^2
\end{align}
Here, $\overline{L}_{\delta}(u)$ consists of all terms which come from one derivative acting on $\delta u$ and the remaining terms acting on $g_{\delta}$. In particular, $\overline{L}_{\delta}(u)$ is linear in $u$ with $O(\delta^2)$ coefficients. $Q_{\delta}(\delta u)$ is at least quadratic in $\delta u$ with bounded coefficients, and $W_{\delta} = F_{\delta}(0) - J_{Y_{\gamma,0}}^{-1}(H(Y_{\gamma}(0), g_0))$ consists of all terms which have no $u$ dependence and hence are $O(\delta)$ and computable in derivatives of $g_{\delta}$. We also note that 
\begin{align*}
J_{Y_{\gamma,0}}:& x^{\mu} \Lambda^{2,\alpha}_0 \to x^{\mu} \Lambda^{0,\alpha}_0 \\
\overline{L}_{\delta} &: x^{\mu} \Lambda^{2,\alpha}_0 \to x^{\mu} \Lambda^{0,\alpha}_0 \\
Q_{\delta}: & x^{\mu} \Lambda^{2,\alpha}_0 \to x^{\mu} \Lambda^{2,\alpha}_0
\end{align*}
are all bounded operators preserving the weight, $x^{\mu}$. Note that the $O(\delta)$ coefficient in equation \eqref{MCExpansionEqn} is given by
\[
C_1 = J_Y(u) + L_2(\dot{g})
\]
in the appendix \S \ref{L2Append}, we show that $L_2 = x^2 \overline{L}_2$ where $\overline{L}_2$ is an order preserving edge operator that maps $2$-tensors to functions. Since $\dot{g} = O(x^{-2})$, this means that $L_2(\dot{g}) = O(1)$. Furthermore, $J_Y$ is an elliptic edge operator for any minimal surface (see \cite{marx2021variations} for computation of the symbol) and hence order preserving. This means that in order for $C_1 = 0$, we must have $u = O(1)$, and so we will choose our weight, $\mu$, to be any $\mu \in (-1,0]$. \nl \nl
We now compute
\begin{align*}
F_{\delta}(u) - F_{\delta}(v) &= -[u - v] - \delta^{-1}J_{Y_{\gamma,0}}^{-1}[\overline{L}_{\delta}(u - v) + Q_{\delta}(\delta u) - Q_{\delta}(\delta v)] + (u - v) \\
&= - \delta^{-1} J_{Y_{\gamma,0}}^{-1}[\overline{L}_{\delta}(u - v)] - \delta^{-1} J_{Y_{\gamma,0}}^{-1}[Q_{\delta}(\delta u) - Q_{\delta}(\delta v)]
\end{align*}
Using that $J_{Y_{\gamma,0}}^{-1}$ is defined and bounded on our weighted spaces (by non-degeneracy), along with our bounds on $\overline{L}_{\delta}$ and $Q_{\delta}$ and the linearity of $\overline{L}$, we have
\begin{align} \label{ContractionEqn}
||F_{\delta}(u) - F_{\delta}(v)|| & \leq K \delta ||u - v||
\end{align}
For $K$ independent of $\delta$. To fully run the fixed point argument, note that $|F_{\delta}(0)| \leq C$ for some $C$ independent of $\delta$ so consider $F_{\delta}: B_{2C}(0) \to B_{2C}(0)$. Equation \eqref{ContractionEqn} shows that our rate of contraction is sufficiently small so there exists a fixed point by the Banach Fixed Point theorem. \qed \nl \nl
\rmk: We will make use of proposition \ref{ExistenceOfSmall} when $Y_{\gamma,0} = HS^{2m^*} \subseteq \H^{n+1}$, which is non-degenerate as is discussed in \cite{alexakis2010renormalized}. We will also apply proposition \ref{ExistenceOfSmall} to small perturbations of $HS^{2m^*}$, which will also be non-degenerate since non-degeneracy is an open condition.
\section{Proof of theorem \ref{ConformalInfinityThm}} \label{ProofSection}
In this section we prove 
\begin{theorem*}
Suppose that $g$ is partially even of order at least $2$, then the renormalized area on all minimal surfaces, $Y^2 \subseteq M^{n+1}$, determines the conformal infinity, $c(g)$, of the ambient manifold.
\end{theorem*}
\noindent We note that the condition of $\omega_1 = 0$ means that renormalized area on $Y^2 \subseteq M^{n+1}$ is well defined and independent of the choice of special bdf. We handle the cases of $n = 2$ and $n \geq 3$ separately.
\subsection{$n = 2$}
\begin{figure}[h!]
	\centering
	\includegraphics[scale=0.4]{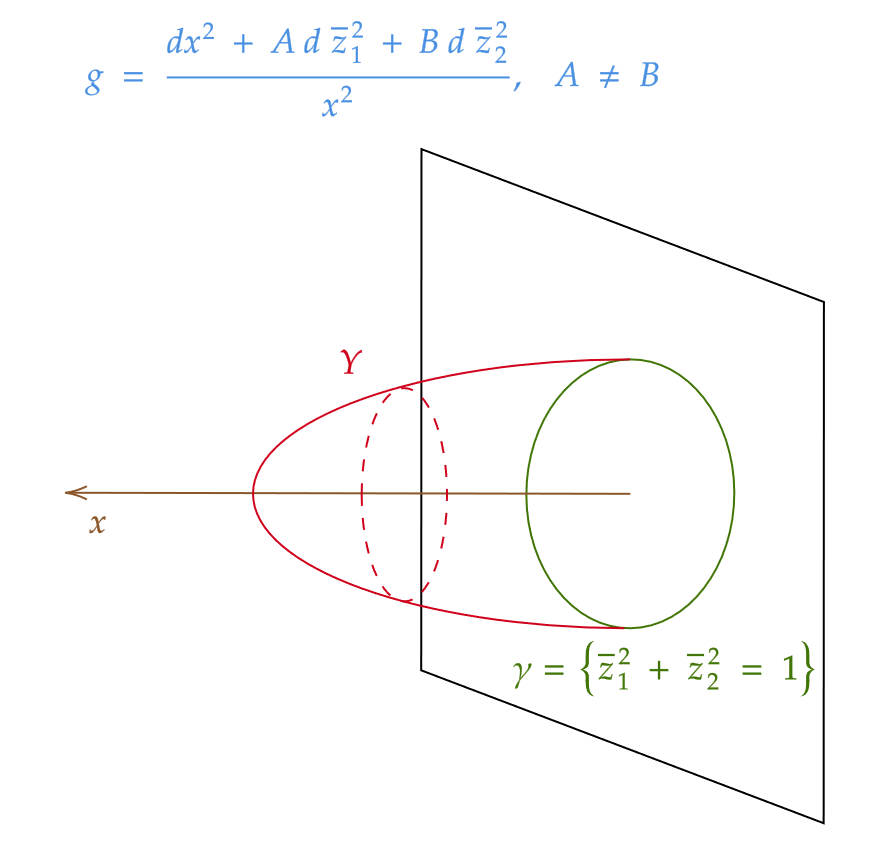}
	\caption{Visualization of the skew minimal surface over a circle in $\overline{z}$ coordinates}
	\label{fig:skewsurface}
\end{figure}

\begin{figure}[h!]
	\centering
	\includegraphics[scale=0.3]{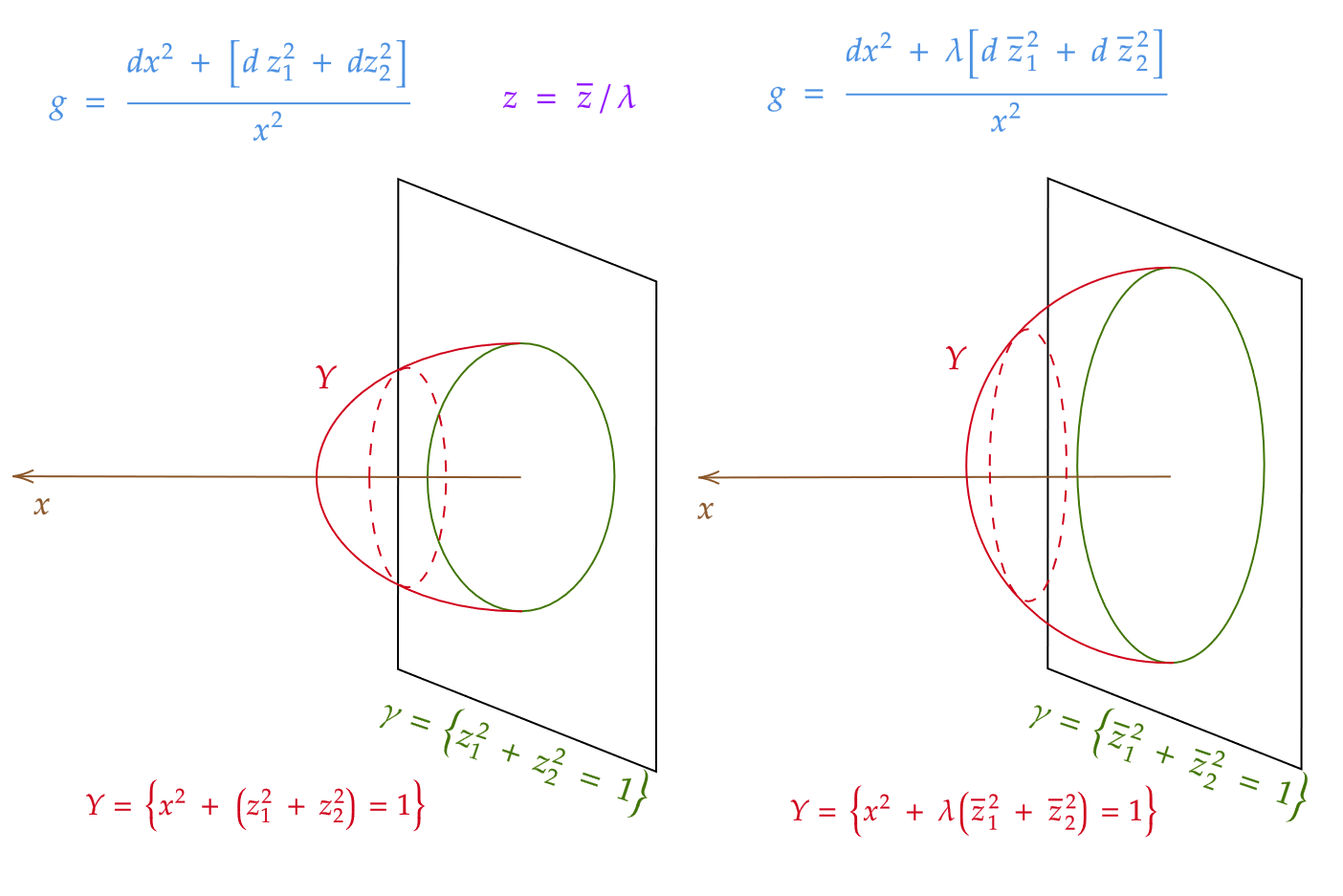}
	\caption{Visualization of hemisphere in different coordinate systems.}
	\label{fig:COfC}
\end{figure}

Given $p \in \partial M$, choose coordinates $y_1, y_2$ for the boundary and let 
\[
\gamma_{\delta}(y) = \{y_1^2 + y_2^2 = \delta^2\}
\]
For each set of coordinates $y = \{y_1, y_2\}$. If we pull back this curve by $F_{\delta}$, we define
\[
F_{\delta}^*(\gamma_{\delta}(y)) = \gamma(\ty) = \{\ty_1^2 + \ty_2^2 = 1\}
\]
And there exists $Y_{\gamma(\ty), \delta}$ which is minimal over $\gamma(\ty)$ with respect to $g_{\delta}$ according to proposition \ref{ExistenceOfSmall}. We then define 
\[
Y_{\gamma_{\delta}(y), \delta}: = F_{\delta}^{-1}(Y_{\gamma(\ty), \delta})
\]
which is minimal with respect to $g$. Thus, we can compute
\begin{align} \label{RAPullback}
\RA(Y_{\gamma_{\delta}(y), \delta}) &= \FPz \int_{Y_{\gamma_{\delta}(y), \delta}} x^z dA_{Y_{\gamma_{\delta}(y), \delta}} \\ \nonumber
&= \FPz \int_{Y_{\gamma(\ty), \delta}} \delta^{z} \tx^z dA_{Y_{\gamma(\ty), \delta}} \\ \nonumber
&= \FPz (1 + z \ln(\delta) + O(z^2)) \int_{Y_{\gamma(y), \delta}} \tx^z dA_{Y_{\gamma(\ty), \delta}} \\ \nonumber
&= \FPz \int_{Y_{\gamma(\ty), \delta}} \tx^z dA_{Y_{\gamma(\ty), \delta}} +  \ln(\delta) \FPz z \int_{Y_{\gamma(y), \delta}} \tx^z dA_{Y_{\gamma(\ty), \delta}} \\ \nonumber
&= A(\ty, \delta) + \ln(\delta) B(y, \delta)
\end{align}
Note that 
\begin{align*}
B(y, \delta) &= \FPz z \int_{Y_{\gamma(y), \delta}} \tx^z dA_{Y_{\gamma(\ty), \delta}}  \\
&= \int_{\gamma(y)} [\sqrt{\det h}]^{1} dA_{\gamma(y)}
\end{align*}
here, we've recalled equation \eqref{AreaFormEqn}
\[
dA_{Y_{\gamma(\ty), \delta}} = \frac{1}{x^m}\sqrt{\det \bh} dx dA_{\gamma(y)}
\]
following the notation and techniques in \cite{marx2021variations}. One can easily check that  
\[
\sqrt{\det h}(s,x) = 1 + O(x^2)
\]
where $s$ is an arclength parameterizing coordinate on $\gamma(y)$ (see equation \eqref{AreaExpansion} and surrounding references, or \cite{graham1999volume}, equation 4.3). Thus
\[
B(y, \delta) = \int_{\gamma(y)} [\sqrt{\det h}]^{1} dA_{\gamma(y)} = 0
\]
%
Thus
\begin{align*}
\RA(Y_{\gamma_{\delta}(y), \delta}) &= A(\ty, \delta) \\
&= \FPz \int_{Y_{\gamma(\ty), \delta}} \tx^z dA_{Y_{\gamma(\ty), \delta}} \\
&= \RA(Y_{\gamma(\ty), \delta}, g_{\delta})
\end{align*}
i.e. our quantity of interest is the same as the renormalized volume of our minimal surface rescaled by $\delta^{-1}$, computed with respect to $g_{\delta}$. In particular, by proposition \ref{ExistenceOfSmall}, we know that $\RA(Y_{\gamma(\ty), \delta}, g_{\delta})$ is continuous with respect to $\delta$, and so we define
\[
A(y, 0) = \lim_{\delta \to 0} A(\ty, \delta) = \RA(Y_{\gamma(\ty), 0}, g_0)
\]
Note that 
\begin{equation} \label{LimitMetric}
g_0 = \frac{d\tx^2 + \omega^0_{ij}(p) d\ty^i d \ty^j}{\tx^2}
\end{equation}
Since $\omega_0(p)$ is positive definite, this is the metric for hyperbolic space after changing coordinates. Thus we can apply equation \eqref{GaussBonnet} so that 
\[
A(y,0) = \RA(Y_{\gamma(\ty), 0}, g_0) =  - 2\pi \chi(Y) - \frac{1}{2} \int_Y |\hat{A}|^2 dA_Y
\]
having noted that the Weyl curvature tensor vanishes in hyperbolic space. Moreover, we know that because $\gamma(\ty)$ is topologically a circle, that $Y_{\gamma(\ty), 0}$ is topologically a disk (see \cite{polthier1994geometric} among many sources). Thus
\[
\RA(Y_{\gamma(\ty), 0}, g_0) = - 2 \pi - \frac{1}{2} \int_Y |\hat{A}|^2 dA_Y
\]
We now aim to detect $\omega_0(p)$ up to a scalar factor by making the above integral vanish. Using corollary \ref{Rigidity}, we know that $\RA(Y_{\gamma(\ty), 0}, g_0) = - 2 \pi $ if and only if $\gamma(\ty)$ is a circle with respect to the boundary metric of $g_0\Big|_{T \partial \H^{n+1}} = \omega_0$. Formally, for a given coordinate system $z = \{z_1, z_2\}$, there exists an orthonormal change of basis to $\ty = \{\tz_1, \tz_2\}$ so that $\gamma(z) = \gamma(\tz)$ and 
\begin{equation} \label{diagonalMetric}
\omega_0 = \omega^0_{ij} dz^i dz^j = \tilde{\omega}^0_{11} (d \tz^1)^2 + \tilde{\omega}^0_{11} (d \tz^2)^2
\end{equation}
See figure \ref{fig:skewsurface} for a visualization of the corresponding minimal surface.
By varying $z_1, z_2$ 
\begin{align*}
\bz_1 &= a z_1 + b z_2 \\
\bz_2 &= c z_1 + d z_2
\end{align*}
there exists some $a,b,c,d$ so that 
\begin{equation} \label{scaleMetric}
	\omega_0 = \lambda [d\bz_1^2 + d\bz_2^2]
\end{equation}
This coordinate system is not unique, as any orthonormal change of basis times a dilation will achieve the same effect. However, in any such coordinate basis, $\{\bz_1, \bz_2\}$, where equation \eqref{scaleMetric} holds, we have
\[
Y_{\gamma(z), 0} = \{x^2 + \frac{\bz_1^2 + \bz_2^2}{\lambda} = 1 \}
\]
is exactly a hemisphere with respect to $g_0$, see figure \ref{fig:COfC}.
As a result,
\[
\RA(Y_{\gamma(z),0}) = - 2\pi
\]
Moreover, if $\bz$ is a coordinate basis which \textit{does not} satisfy equation \eqref{scaleMetric}, then after an orthonormal transformation (which again, preserve $\gamma(\by)$), $\omega_0$ can be written according to equation \eqref{diagonalMetric} but $Y_{\gamma(\bz),0}$ will not be a geodesic hemisphere with respect to $g_0$ and so $\RA(Y_{\gamma(\bz),0}) < - 2 \pi$. This means that we can detect a coordinate system for which $\omega_0(p)$ satisfies equation \eqref{scaleMetric} by measuring when $\RA(Y_{\gamma(z),0}) = - 2\pi$ exactly. \nl \nl
Once such a coordinate system, $\ty$, is found in the rescaled setting, we know it comes from blowing up coordinates on $\partial M$ near $p$ via the map $F_{\delta}$ from equation \eqref{DeltaDiffeo}. Thus, there exists coordinates $y$ on a neighborhood of $p \in \partial M$ which are orthonormal with respect to $\omega_0(p)$, i.e. $\omega_0(p) = \lambda g_{\text{euc},2}$. \nl \nl
Given $q \in \partial M$ near $p$, we can find another coordinate system $z(q) = \{z_1(q), z_2(q)\}$ so that $\omega_0(q) = \lambda(q) g_2(z(q))$, i.e. $\omega_0$ is a multiple of the identity with respect to $\p_{z_1(q)}, \p_{z_2(q)}$. Let $z(p) =  \{y_1, y_2\}$. We can define an auxiliary metric 
\[
\overline{\omega}_0(q) = g_{\text{euc},2}(z(q))
\]
i.e. our metric is that which is the identity in the $\{\p_{z_1(q)}, \p_{z_2(q)}\}$ basis at each point. The issue is that this might not be continuous, so we define
\[
\tilde{\omega}_0(q) := s(q) g_2(z(q))
\]
for some scalar valued function, $s(q)$, so that $\tilde{\omega}_0(q)$ is a \textit{continuous} metric on the boundary near $p$. Evidently, such an $s(q)$ exists since $s(q) = \lambda(q)$ would give $\tilde{\omega}_0(q) = \omega_0(q)$ which is smooth by assumption. However, $s(q) = f(q) \cdot \lambda(q)$ for any positive continuous $f(q)$ also works, reflecting that we can only hope to construct an element of the conformal infinity. \qed

\subsection{$n \geq 3$}
The proof is identical up to equation \eqref{diagonalMetric}. From then on, our construction of minimal surfaces allows us to analyze $\omega_0$ restricted to $2$-dimensional subspaces. In particular, since our construction of $Y_{\gamma(\ty), 0}$ has $\gamma(\ty) \subseteq \R^2_{\ty_1, \ty_2} \subseteq \R^n$ we know that $Y_{\gamma(\ty), 0} \subseteq \R^2_{\ty_1, \ty_2} \times \R^+ \subseteq \R^n \times \R^+$. Thus it suffices to consider the behavior of 
\[
\omega_0 \Big|_{\R^2_{\ty_1, \ty_2}} = \omega^0_{11} (d \ty^1)^2 + \omega^0_{12} d \ty^1 d \ty^2 + \omega^0_{22} (d \ty^2)^2
\]
As in the $n = 2$ case, we can vary $\ty_1, \ty_2$ within their own span to detect when 
\begin{equation} \label{y1y2Basis}
\omega_0 \Big|_{\R^2_{\ty_1, \ty_2}} = \lambda [(d \ty^1)^2  + (d \ty^2)^2]
\end{equation}
Now choose a coordinate $\ty_3$ which is linearly independent from $\ty_1, \ty_2$, and consider the pair of coordinates $\{\ty_1, \ty_3\}$ with 
\[
\omega_0 \Big|_{\R^2_{\ty_1, \ty_3}} = \omega^0_{11} (d \ty^1)^2 + \omega^0_{13} d \ty^1 d \ty^3 +  \omega^0_{33} (d \ty^3)^2
\]
We now \textit{fix} $\ty_1$ and \textit{vary} $\ty_3 \rightarrow  A \ty_1 + B \ty_3 = \ty_3'$ until
\begin{equation} \label{scaleMetric2}
\omega_0 \Big|_{\R^2_{\ty_1, \ty_3}} = \overline{\lambda} [(d \ty^1)^2  + (d \ty^3)^2]
\end{equation}
\noindent Such a change of coordinates exists by setting, e.g.
\begin{align*}
B &= -A \frac{\omega_{13}}{\omega_{11}} \\
A &= \sqrt{\frac{\omega_{11} \omega_{33} - \omega_{13}^2}{\omega_{13}}}
\end{align*}
This presumes $\omega_{13} \neq 0$ - if $\omega_{13} = 0$, we can just set $\ty_3' = \sqrt{\omega_{11}^0 / \omega_{33}^0} \ty_3$, and $S^1$ with respect to $\ty_1$ and $\ty_3'$ will genuinely be a circle, as detected by $\RA(Y_{\gamma(\ty_1, \ty_3'),0}) = - 2\pi$. \nl \nl
With this, we have a basis $\{\ty_1, \ty_2, \ty_3\}$ so that 
\[
\omega_0 \Big|_{\R^3_{\ty_1, \ty_2, \ty_3}} = \lambda [(d \ty^1)^2 + (d \ty^2)^2 + (d \ty^3)^2] + A d \ty^2 d \ty^3 \sim \begin{pmatrix}
	\lambda & 0 & 0 \\
	0 & \lambda & A \\
	0 & A & \lambda
	\end{pmatrix}
\]
Here, we've noted that $\overline{\lambda} = \lambda$ from equations \eqref{y1y2Basis} and \eqref{scaleMetric2} since $\ty_1$ was determined independent of $\ty_3$. In particular, $\omega_0(\p_{\ty_1}, \p_{\ty_2}) = \omega_0(\p_{\ty_1}, \p_{\ty_3}) = 0$, so we can now repeat the same construction with $\ty_2, \ty_3$, keeping $\ty_2$ fixed and modifying $\ty_3$ until $\RA(Y_{\gamma(\ty_2, \ty_3),0}) = -2\pi$. In particular, since $\ty_1, \ty_2$ are fixed in this construction, we'll have  
\[
\omega_0 \Big|_{\R^3_{\ty_1, \ty_2, \ty_3}} = \lambda [(d \ty^1)^2 + (d \ty^2)^2 + (d \ty^3)^2] 
\]
This construction generalizes by induction: given 
\begin{equation} \label{InductiveHypothesis}
\omega_0 \Big|_{\R^{k-1}_{\ty_1, \ty_2, \dots, \ty_{k-1}}} = \lambda [(d \ty^1)^2 + (d \ty^2)^2 + \dots + (d \ty^{k-1})^2] 
\end{equation}
We can consider the pairs of coordinates $(\ty_j, \ty_k)$. For each $j \leq k - 1$, replace $\ty_k \rightarrow A_j \ty_j + B_j \ty_k$ so that 
\[
\RA(Y_{\gamma(\ty_j, \ty_k), 0}) = - 2\pi
\]
After doing this for all $1 \leq j \leq k-1$, we will have achieved
\[
\omega_0 \Big|_{\R^k_{\ty_1, \ty_2, \dots, \ty_{k}} } = \lambda [(d \ty^1)^2 + (d \ty^2)^2 + \dots + (d \ty^{k})^2]
\]
repeating the same construction to find a smooth conformal representative as in the $n = 2$ case, we complete the proof. \qed

\section{Proof of theorem \ref{ExpansionThm}}
The goal of this section is to prove theorem \ref{ExpansionThm}
\begin{theorem*}
Let $(M^n, g)$, $g$ a partially even metric to order $2m^*$, and $n > 2m^*$. Suppose $\omega_0$ is known in equation \eqref{AHExpansion}. Then the renormalized volume on all $2m^*$-dimensional minimal submanifolds, $Y^m$, determines the asymptotic expansion in equation \eqref{AHExpansion} to arbitrary order.
\end{theorem*}
\noindent \Pf There exists a bdf $x: \overline{M}^{n} \to \R^{\geq 0}$ such that $g$ takes the explicit form of equation \eqref{AHExpansion}
\[
g = \frac{dx^2 + \omega_0(y) + x^2 \omega_2(y) + O(x^3)}{x^2} 
\]
and $\omega_0$ is a representative of the conformal class, $c(g)$. We first handle the codimension $1$ case, i.e. knowledge of $\RV(Y^{2m^*})$ for $Y^{2m^*} \subseteq M^{2m^*+1}$
\subsection{Codimension $1$, Base Case}
\begin{figure}[h!]
	\centering
	\includegraphics[scale=0.3]{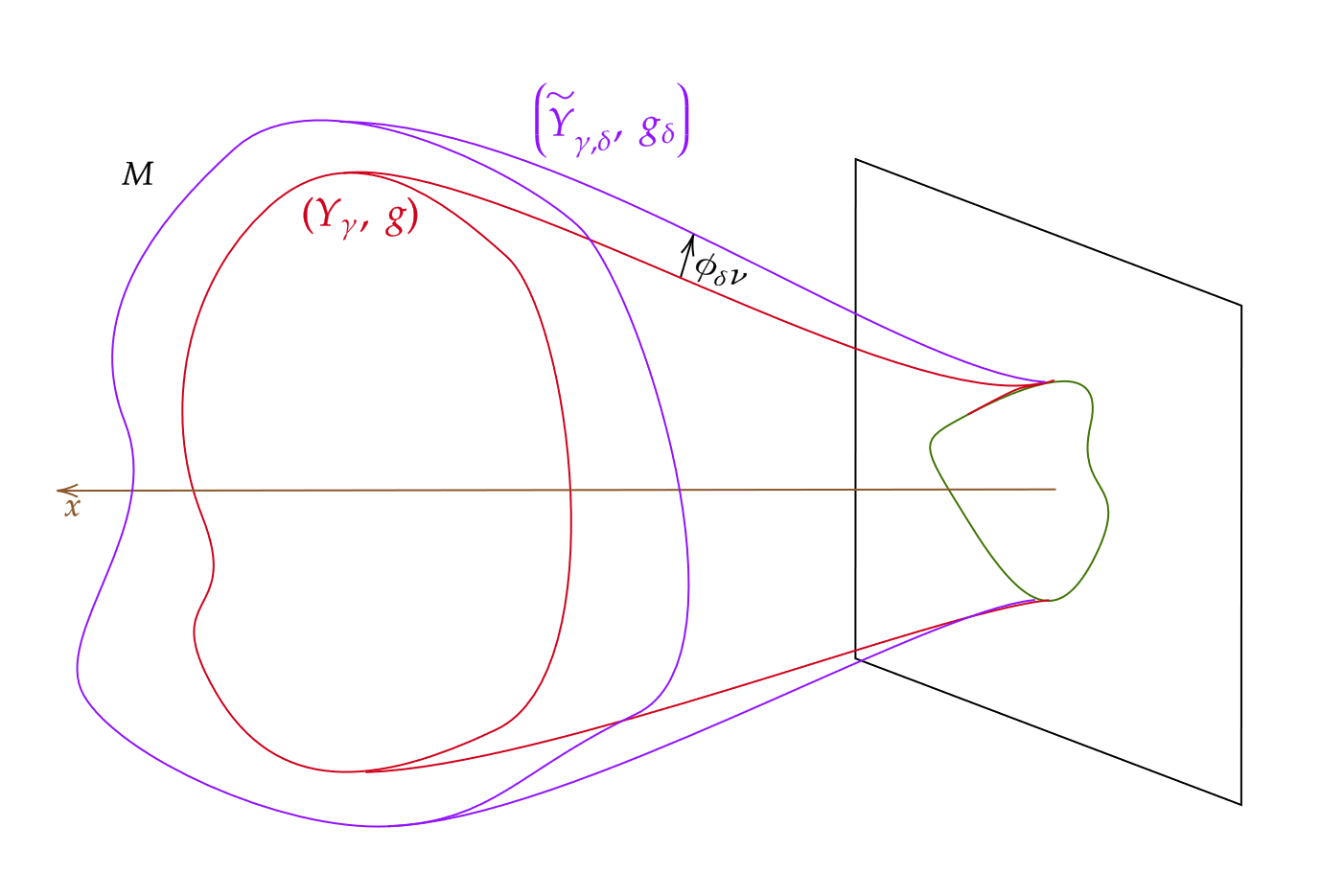}
	\caption{Visualization of $\tY_{\gamma, \delta}$}
	\label{fig:ydeltaminimal}
\end{figure}

As in the proof of theorem \ref{ConformalInfinityThm}, for any point $p \in \partial M$, we can find an open neighborhood $U$ and coordinate system $\{x, y_1, y_2, \dots, y_{2m^*}\}$ so that $\omega_0(p) = [dy_1^2 + \dots dy_{2m^*}^2]$ and $x$ is a special bdf. As in the proof of theorem \ref{ConformalInfinityThm} we perform the following construction: \nl \nl
Let $\gamma^{2m^*-1} \subseteq \R^{2m^*}$ be a closed smooth submanifold with diameter $D_{\gamma}$. There exists a minimal submanifold, $Y_{\gamma}$, with asymptotic boundary $\gamma^{2m^*-1}$ with respect to the standard hyperbolic metric 
\[
g_{\H^{2m+1}} = \frac{d\tx^2 + [d\ty_1^2 + \dots + d \ty_{2m^*}^2]}{\tx^2}
\]
which is topologically compact in $\R^{\geq 0} \times \R^{2m^*}$. Furthermore, assume that $Y_{\gamma}$ is unique and non-degenerate. Define
\[
\gamma(y, \delta) = \{(y_1, \dots, y_{2m^*}) = \delta (\ty_1(p), \dots, \ty_{2m^*}(p)) \; | \; p \in \gamma \subseteq \R^{2m^*}\} \subseteq \partial M
\]
And $Y_{\gamma, \delta} \subseteq M$ to be a minimal surface with asymptotic boundary $\gamma(y, \delta)$. This exists by proposition \ref{ExistenceOfSmall}, and by the maximum principle, any such minimal surface must be contained in a hemisphere of radius $2\delta R$ for some $R > 0$ independent of $\delta$. Note that the hemisphere of (euclidean) radius $2\delta R$ is minimal with respect to $g_{\H^{2m^*+1}}$ but not minimal with respect to $g_{\delta}$. However, we can generate a barrier by considering hemispheres centered at $x = \eta \delta R$, call this $HS^n(\delta R, \eta \delta R)$ with $\eta < 1$ fixed. For $\delta$ sufficiently small, we have 
\[
H_{HS^n(R,\eta R)} \geq \eta + O(\delta)
\]
which will be strictly positive as $\delta \to 0$ and act as a barrier. This is identical to shifted hemispheres constructed in the envelope argument outlined in \S 3.4.2 of \cite{marx2021variations}. \nl 
\indent Having established boundedness of $Y_{\gamma, \delta}$, we can now pull back $Y_{\gamma, \delta}$ to $\R^{\geq 0} \times \R^{2m^*}$ via the dilation map, $F_{\delta}$, from equation \eqref{DeltaDiffeo}. In particular $\tY_{\gamma, \delta}: = F_{\delta}^{-1}(Y_{\gamma, \delta})$ is a minimal surface with respect to the metric $g_{\delta} = F_{\delta}^*(g)$. Moreover, for some $\kappa > 0$,
\[
\tY_{\gamma, \delta} \subseteq U = [0,\delta^{-1} \kappa) \times B_{\kappa \delta^{-1}}^{2m^*}(0) \subseteq \R^{\geq 0} \times \R^{2m^*}
\]
so we can think of $\tY_{\gamma, \delta}$ and $Y_{\gamma}$ lying in the same chart but with different metrics. By proposition \ref{ExistenceOfSmall}, we know that $\tY_{\gamma, \delta}$ exists and is unique and is asymptotic to $\gamma$ (see figure \ref{fig:ydeltaminimal})

We first show that we can recover $\omega_k$ for $k = 2$. Recall that $x = F_{\delta}(\tx) = \delta \tx$. We compute
\begin{align*}
\RV(Y_{\gamma, \delta}, g) &= \FPz \int_{Y_{\gamma, \delta}} x^z dA_{Y_{\gamma, \delta}} \\
&= \FPz \int_{\tY_{\gamma, \delta}} \delta^z \tx^z dA_{\tY_{\gamma, \delta}} \\
&= \FPz \int_{\tY_{\gamma, \delta}} \tx^z dA_{\tY_{\gamma, \delta}} + \ln(\delta) \FPz z \int_{\tY_{\gamma, \delta}} \tx^z dA_{\tY_{\gamma, \delta}}
\end{align*}
having expanded $\delta^z = 1 + z \ln(\delta) + O(z^2)$ and used that the $O(z^2)$ terms vanish under finite part evaluation as described in section \S \ref{mechanics}. Furthermore
\begin{align*}
\FPz z \int_{\tY_{\gamma, \delta}} \tx^z dA_{\tY_{\gamma, \delta}} &= \int_{\gamma} [\sqrt{\det \bh}]^{(2m^*-1)} dA_{\gamma}
\end{align*}
where $\bh$ denotes the induced compactified metric, $\bg_{\delta} = \tx^2 g_{\delta}$ on $\tY_{\gamma, \delta}$. Using that $\sqrt{\det \bh}$ is even up to order $2m^*$ (see equation \eqref{AreaFormEqn}), we see that
\begin{align*}
\FPz z \int_{\tY_{\gamma, \delta}} \tx^z dA_{\tY_{\gamma, \delta}} &= 0
\end{align*}
With this, we have
\begin{align*} 
\RV(Y_{\gamma, \delta}, g) & = \FPz \int_{\tY_{\gamma, \delta}} \tx^z dA_{\tY_{\gamma, \delta}} \\
&= \RV(\tY_{\gamma, \delta}, g_{\delta})
\end{align*}
here, we've noted that $\tx$ is a special bdf for $g_{\delta}$. With this, we can compute
\begin{align*}
\frac{d^2}{d \delta^2} \RA(Y_{\gamma}, g_{\delta}) \Big|_{\delta = 0} &= \frac{d^2}{d \delta^2} \FPz \int_{Y_{\gamma}} \tx^z dA_{Y_{\gamma}} 	
\end{align*}
We now pullback to $Y_{\gamma}$, our minimal surface in hyperbolic space. Let 
\begin{align} \label{PDeltaMap}
P_{\delta}&: Y_{\gamma} \rightarrow \tY_{\gamma, \delta} \\ \nonumber
P_{\delta}(p) &= \exp_{Y_{\gamma}, g_{\H^{2m^*+1}} }(\phi_{\delta}(p) \nu(p))
\end{align}
Then
\begin{align} \label{SecondDeltaDeriv}
\frac{d^2}{d \delta^2} \RA(\tY_{\gamma, \delta}, g_{\delta}) \Big|_{\delta = 0} &= \frac{d^2}{d \delta^2} \FPz \int_{Y_{\gamma}} P_{\delta}^*(\tx)^z P_{\delta}^*(dA_{\tY_{\gamma, \delta}}) \\ \nonumber
&= \FPz \int_{Y_{\gamma}} z \tx^{z-1} \dot{\tx}_{\delta \delta} + z(z-1) \tx^{z-2} \dot{\tx}_{\delta}^2 + 2 \tx^{z-1} \dot{\tx}_{\delta} \frac{d}{d \delta} P_{\delta}^*(dA_{\tY_{\gamma, \delta}}) +  \tx^z \frac{d^2}{d \delta^2} P_{\delta}^*(dA_{\tY_{\gamma, \delta}}) \\ \nonumber
&= A + B + C + D
\end{align}
Here, $\dot{\tx}_{\delta\delta}$ and $\dot{\tx}_{\delta}$ records how the $x$ component of a point on $\tY_{\gamma, \delta}$ moves to first and second order as $\delta \to 0$. Given that 
\[
\mathcal{H}(\tY_{\gamma, \delta}, g_{\delta}) = 0
\]
we can differentiate this and compute $\dot{\tx}_{\delta}$ from $\dot{\phi} = \frac{d}{d \delta} \phi_{\delta} \Big|_{\delta =0}$:
\begin{align} \label{MCVaryingDelta}
\mathcal{H}(\tY_{\gamma, \delta}), g_{\delta}) &= 0 \\ \nonumber
\implies J_Y(\dot{\phi}) + L_2(\dot{g}) &= 0 \\ \label{DotPhiEq}
\implies \dot{\phi} &= - J_Y^{-1}(L_2(\dot{g}))
\end{align}
Here, $L_2$ is the linearization coming from fixing $Y = Y_{\gamma}$ and varying the induced metric, $g_{\delta} \to g_0$ on it. We also assume that $Y_{\gamma}$ is non-degenerate with respect to the hyperbolic metric, so that $J_Y^{-1}$ is well defined. We compute in the appendix \ref{L2Append} that 
\begin{equation} \label{L2Form}
L_2(\dot{g}) = - \langle \dot{g}, A \rangle + \div_Y(\dot{g}(\nu, \cdot)) - \frac{1}{2} \nu(\tr_Y(\dot{g}))
\end{equation}
Equation \eqref{MCVaryingDelta} similarly gives
\begin{align} \label{SecondVarEqn} 
\mathcal{H}(P_{\delta}(Y_{\gamma}), g_{\delta}) &= 0 \\ \nonumber
\implies J_Y(\ddot{\phi}) + L_2(\ddot{g}) + DJ_Y(\dot{g})(\dot{\phi}) + DL_2(\dot{\phi})(\dot{g}) &= 0 
\end{align}
where $DJ_Y(\dot{g})$ means the derivative of the Jacobi operator on $Y$ as the metric on $Y$ changes, and $DL_2(\dot{\phi})$ means the change in $L_2$ as we infinitesimally vary the surface $Y$ by $\dot{\phi}$. \nl \nl
In general, we will only need the existence of $\dot{\phi}^{(k)}$ which can be computed in terms of $k$ derivatives of the metric $g_{\delta}$. See \cite{marx2021variations} for more information on the regularity theory of these variations. We now analyze each term in equation \eqref{SecondDeltaDeriv}.
%
%
%
%
\subsubsection{A} \label{ASection}
We compute
\begin{align*}
A &= \FPz \int_{Y_{\gamma}} z \tx^{z-1} \dot{\tx}_{\delta \delta} dA_{Y_{\gamma}} \\
&= \int_{\gamma} [\dot{\tx}_{\delta \delta} \bJac_Y]^{(m)} dA_{\gamma}
\end{align*}
where $\bJac_Y = \sqrt{\det \bh}$ for $\bh = \bg \Big|_{TY}$. 
\subsubsection{B}
Since $\dot{\tx}_{\delta}$ just depends on $\dot{\phi}$, which is dependent on $\dot{g}$, we know that $B = F(\omega_0, \nabla \omega_0)$, i.e. it is known by assumption. 
\subsubsection{C}
We compute
\begin{align*}
\frac{d}{d \delta} P_{\delta}^*(dA_{Y_{\gamma, \delta}}) &= [\dot{\phi} H_Y + \tr_Y(\dot{g})] dA_Y \\
&= \tr_Y(\dot{g}) dA_Y
\end{align*}
so using equation \eqref{DotPhiEq} and that $Y$ is minimal with respect to the hyperbolic metric, $g_0$, we see that $C = F(\omega_0, \nabla \omega_0)$.
\subsubsection{D} \label{DTermAnalysis}
We compute
\begin{align*}
D &= \FPz \int_{Y_{\gamma}} \tx^z \frac{d^2}{d \delta^2} P_{\delta}^*(dA_{Y_{\gamma, \delta}})
\end{align*}
Explicitly, we know that
\begin{align*}
\frac{d}{d \delta} P_{\delta}^*(dA_{Y_{\gamma, \delta}}) \Big|_{\delta = \delta_0} &= \left[- \langle \dot{P}_{\delta_0}, H_{Y_{\gamma, \delta_0}} \rangle + \div_{Y_{\gamma, \delta_0}}(\dot{P}_{\delta_0}^T) + \tr_{Y_{\gamma, \delta_0}}(\dot{g}_{\delta_0})\right] dA_{Y_{\gamma, \delta_0}}\\
\implies \frac{d^2}{d \delta^2} P_{\delta}^*(dA_{Y_{\gamma, \delta}}) \Big|_{\delta = 0} &= \Big(\left[ - \langle A, \dot{P}\rangle^2 + |\n_{Y_{\gamma}}^N \dot{P}|^2 - \tr_{Y_{\gamma}}(R(\cdot, \dot{P}, \cdot, \dot{P})) + \div_{Y_{\gamma}}(\ddot{P}) \right] \\
& + [\tr_{Y_{\gamma}}(\ddot{g})]  + K(\omega_0, \nabla \omega_0, \nabla^2 \omega_0) \Big) dA_{Y_{\gamma}} \\
&= [\div_{Y_{\gamma}}(\ddot{P}) + \tr_{Y_{\gamma}}(\ddot{g})] dA_{Y_{\gamma}} + K(\omega_0, \n \omega_0, \n^2 \omega_0)
\end{align*}
This can be seen as follows: the pull back of the area form is a function of what $Y_{\gamma(\ty), \delta}$ is as a surface and the metric, $g_{\delta}$, which is restricted to that surface. Since we only care the summands terms which record $\ddot{g} = \omega_2 + K(\omega_0)$ (as opposed to those that depend on just $\dot{g}$ and $g_0$, which are implicitly functions of just $\omega_0$), we can isolate the last two terms and let the remainder be a function of $\omega_0$. Note that 
\begin{align*}
\int_{Y_{\gamma}} \tx^z \div_{Y_{\gamma}}(\ddot{P}) &= - \int_{Y_{\gamma}} z \tx^{z-1} d\tx(\ddot{P}) \\
&= -\int_{Y_{\gamma}} z \tx^{z-1} \dot{\tx}_{\delta \delta} \\
\end{align*}
for $\Re(z) \gg 0$. Now taking the finite part, we have
\begin{align*}
\FPz \int_{Y_{\gamma}} \tx^z  \div_{Y_{\gamma}}(\ddot{P}) &= - \int_{\gamma} [ \dot{\tx}_{\delta \delta} \bJac_Y]^{(m)}
\end{align*}
and so this term cancels with term $A$ in section \ref{ASection}. 
\subsubsection{In Sum}
In sum, we have for $Y_{\gamma}$ a minimal surface over an arbitrary submanifold $\gamma^{2m-1} \subseteq \R^{2m}$:
\begin{align} \label{SecondDerivEqn}
\frac{d^2}{d \delta^2} \RA(Y_{\gamma, \delta}, g_{\delta}) \Big|_{\delta = 0} &= \frac{d^2}{d \delta^2} \FPz \int_{Y_{\gamma}} P_{\delta}^*(\tx)^z P_{\delta}^*(dA_{Y_{\gamma, \delta}}) \\  \nonumber
&= \FPz \int_{Y_{\gamma}} z \tx^{z-1} \dot{\tx}_{\delta \delta} + z(z-1) \tx^{z-2} \dot{\tx}_{\delta}^2 + 2 \tx^{z-1} \dot{\tx}_{\delta} \frac{d}{d \delta} P_{\delta}^*(dA_{Y_{\gamma, \delta}}) +  \tx^z \frac{d^2}{d \delta^2} P_{\delta}^*(dA_{Y_{\gamma, \delta}}) \\ \nonumber
&= \FPz \int_{Y_{\gamma}} \tx^z \tr_{Y_{\gamma}}(\omega_2) dA_{Y_{\gamma}} + K(\omega_0)
\end{align}
We now note that the above recovers the trace of $\omega_2$ at $p$
\begin{lemma} \label{TraceRecoveryLemma}
Equation \eqref{SecondDerivEqn} recovers the trace of $\omega_2$ at $p$:	
\end{lemma}
\noindent \Pf We set $\gamma = S^{2m^*-1}$ and $Y = HS^{2m^*}$ in equation \eqref{SecondDerivEqn}. Then using the computations from \S \ref{HemisphereParamSection} we have
\begin{align*}
\tr_{Y_{\gamma(\ty),0}}(\omega_2) &= \tr_{HS^{2m^*}}(\omega_2) \\
&= h^{\tx\tx} \omega_2(v_{\tx}, v_{\tx}) + h^{ij} \omega_2(v_i, v_j) \\
&= \tx^4 \omega_2(F_{2m^*-1}(\theta), F_{2m^*-1}(\theta)) + \tx^2 \tr_{S^{2m^*}-1}(\omega_2)(\theta) \\
dA_{HS^{2m^*}} &= \frac{(1-\tx^2)^{m^*-1}}{\tx^{2m^*}}
\end{align*}
And so 
\begin{align} \label{TracePrelimEqn}
\frac{d^2}{d \delta^2} \RA(Y_{S^{m-1}, \delta}, g_{\delta}) \Big|_{\delta = 0} &= \FPz \int_{HS^{2m^*}} \tx^z \tr_{HS^{2m^*}}(\omega_2) dA_{HS^{2m^*}} \\ \nonumber
&= [\omega_{2,11} + \dots + \omega_{2,(2m^*, 2m^*)}] \left( c_1 I_1 + c_2 I_2 \right)
\end{align}
where 
\begin{align*}
c_1 [\omega_{2,11} + \dots + \omega_{2,(2m^*, 2m^*)}] &= \int_{S^{2m^*-1}}\omega_2(F_{2m^*-1}(\theta), F_{2m^*-1}(\theta)) dA_{S^{2m^*-1}} \\
c_2 [\omega_{2,11} + \dots + \omega_{2,(2m^*, 2m^*)}] &= \int_{S^{2m^*-1}} \tr_{S^{2m^*}-1}(\omega_2)(\theta) dA_{S^{2m^*-1}} \\
I_1 &= \FPz \int_0^1 \tx^{z+4-2m^*}(1-\tx^2)^{m^*-1} \\
I_2 &= \FPz \int_0^1 \tx^{z+2-2m^*}(1-\tx^2)^{m^*-1}
\end{align*}
We parameterize 
\[
S^{2m^*-1} = \{(y_2, \dots, y_{2m^*}) = (\cos \phi, \sin \phi F_{m-2}(\theta_{m-2})) \} 
\]
and compute
\begin{align*}
c_1 &= \frac{1}{2m^*-1} \left[\int_{\phi = 0}^{\pi} (\sin \phi)^{2m^*}d\phi\right] \text{Vol}(S^{2m^*-2}) \\
c_2 &= (2m^*-1) c_1
\end{align*}
We further compute $I_1, I_2$ via integration by parts:
\begin{align*}
I_1 &= \FPz \int_0^1 \tx^{z+4-2m^*}(1-\tx^2)^{m^*-1} \\
&= \FPz \frac{2(m^*-1)}{z+5-2m^*} \int_0^1 \tx^{z+6-2m^*}(1-\tx^2)^{m^*-2} \\
&= \FPz \frac{2^{m^*-1}(m^*-1)!}{(z+5-2m^*)(z+7-2m^*)\cdots(z-1) (z + 1)} \int \tx^{z+2} \\
&= \frac{2^{m^*-1} (m^*-1)!}{(5-2m^*)(7-2m^*) \cdots (-1) (1) (3)} 
\end{align*}
Having noted that after integrating by parts $m^*-2$ times, there was no pole at $z = 0$ due to parity. Similarly
\begin{align*}
I_2 &= \FPz \int_0^1 \tx^{z+2-2m^*}(1-\tx^2)^{m^*-1} \\
&= \FPz \frac{2(m^*-1)}{(z+3-2m^*)} \int_0^1 \tx^{z+4-2m^*}(1-\tx^2)^{m^*-2} \\
&= \FPz \frac{2^{m^*-1} (m^* - 1)!}{(z + 3 - 2m^*) (z + 5 - 2m^*) \cdots (z-1)} \int_0^1 \tx^z \\
&= \frac{2^{m^* - 1} (m^* - 1)!}{(3 - 2m^*)(5 - 2m^*) \cdots (-1)} \cdot 1 \\
&= \frac{3}{2m^*-3} I_2
\end{align*}
%
%
Together, this gives
\begin{align*}
c_1 I_1 &= \frac{c_2}{2m^*-1} \cdot \left[ (-1) \frac{3-2m^*}{3} \right] I_2 \\
&= c_2 I_2 \frac{2m^* - 3}{(2m^*-1) 3} \\
\implies c_1 I_1 + c_2 I_2 &= c_2 I_2 \left[1 + \frac{(2m^*-3)}{3(2m^*-1)} \right] \\
&\neq 0
\end{align*}
This computation along with equation \eqref{TracePrelimEqn} finishes the proof. \qed 
\subsubsection{Varying $\gamma$}
\begin{figure}[h!]
	\centering
	\includegraphics[scale=0.3]{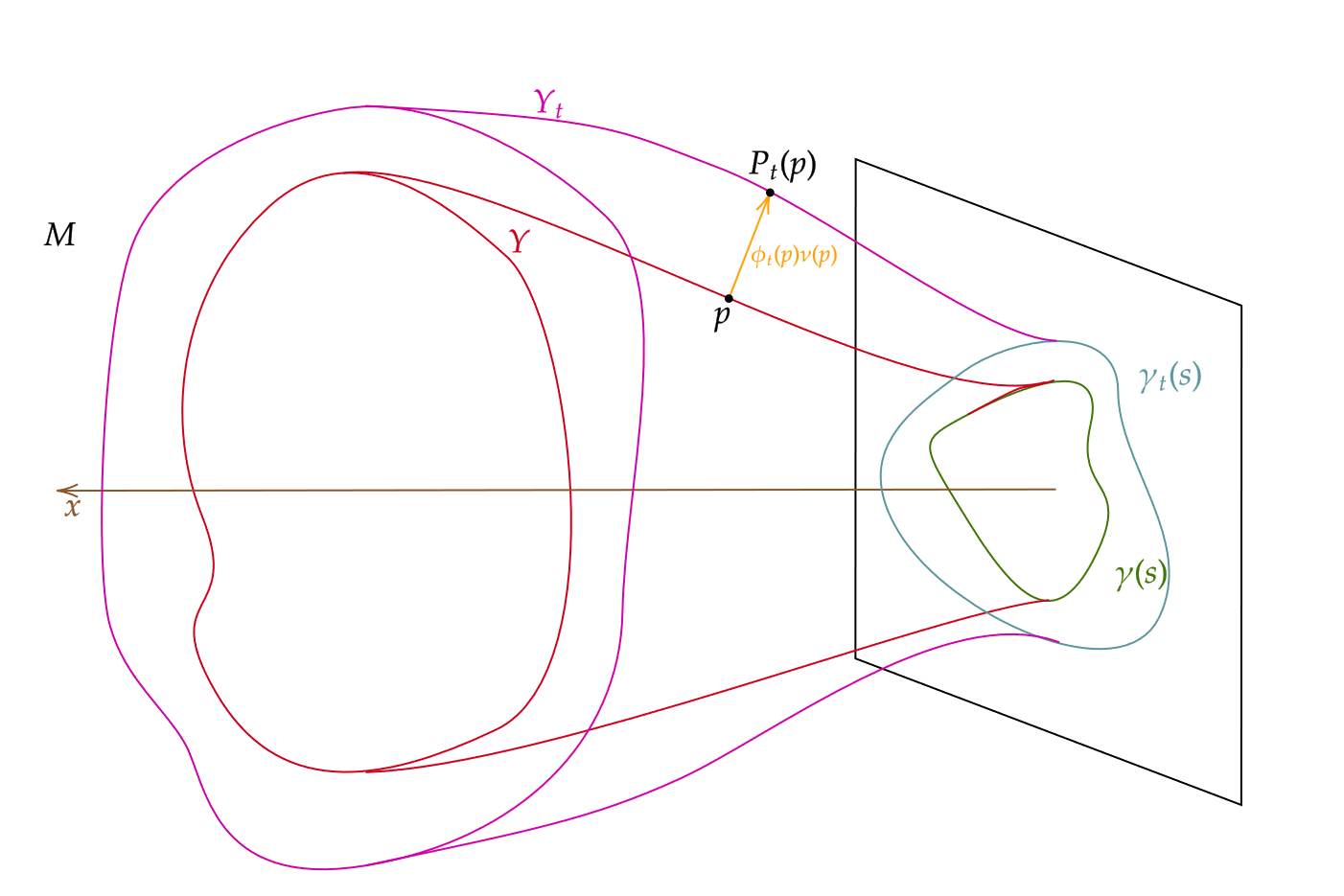}
	\caption{Visualization of $Y_t$ over $Y$ and $\gamma_t$ over $\gamma$.}
	\label{fig:ytovery}
\end{figure}

Having recovered the trace, we now set $\gamma = \gamma_t$ to be a perturbation of $S^{m-1} = \{\ty_1^2 + \ty_2^2 + \dots + \ty_{2m}^2 = 1\}$ in equation \eqref{SecondDerivEqn}. Let $Y_{\gamma(t)} \subseteq \H^{2m^*+1}$ be the minimal hypersurface asymptotic to $\gamma(t)$. This follows by the non-degeneracy of $Y_{\gamma}$ (geodesic hemispherical copy of $\H^{2m} \subseteq \H^{2m+1}$) which implies that $Y_{\gamma(t)}$ exists and is unique. Define a graphical parameterization of $Y_{\gamma(t)}$ over $Y_{\gamma}$:
\begin{align*}
P_t &: Y_{\gamma} \rightarrow Y_{\gamma_t} \\
P_t(p) &= \exp_{Y_{\gamma}}(\phi_t(p) \nu(p))	
\end{align*}
in analogy with equation \eqref{PDeltaMap}, see figure \ref{fig:ytovery}.
We can then pull back our variations (with respect to $\delta$) of renormalized area to $Y_{\gamma}$ and differentiate in $t$
\begin{align*}
\frac{d^2}{d \delta^2} \RA(Y_{\gamma_t, \delta}, g_{\delta}) \Big|_{\delta = 0}  &= \FPz \int_{Y_{\gamma_t}} \tx^z \tr_{Y_{\gamma_t}}(\omega_2) dA_{Y_{\gamma_t}} + K(\omega_0, t) \\
&= \FPz \int_{Y_{\gamma}} P_t^*(\tx^z) P_t^*[\tr_{Y_{\gamma_t}}(\omega_2)] P_t^*[dA_{Y_{\gamma_t}}] + K(\omega_0, t) \\
\implies \frac{d}{dt} \frac{d^2}{d \delta^2} \RA(Y_{\gamma_t, \delta}, g_{\delta}) \Big|_{\delta = 0} \Big|_{t = 0}&= \FPz \int_{Y_{\gamma}}z \tx^{z-1} \dot{x}_{\phi} \tr_{Y_{\gamma}}(\omega_2) dA_{Y_{\gamma}} + \tx^z \left(\frac{d}{dt} P_t^*[\tr_{Y_{\gamma_t}}(\omega_2)] \right) dA_{Y_{\gamma_0}} \\
&= A + B
\end{align*}
having noted that differentiating the area form with respect to $t$ yields $0$, since we are only working with minimal surfaces and our variation is normal.  Here $\dot{x}_{\phi} = \dot{\phi} dx(\nu)$ denotes the change in $x$ along $P_t(p)$. We compute this variation in parts.
\subsubsection{$A$}
We compute
\begin{align*}
A& = \FPz \int_{Y_{\gamma}}z \tx^{z-1} \dot{\tx}_{\phi} \tr_{Y_{\gamma}}(\omega_2) dA_{Y_{\gamma}} \\
&= \int_{\gamma} [ \dot{\tx}_{\phi} \tr_{Y_{\gamma}}(\omega_2) \bJac_{Y_{\gamma}}]^{(2m^*)}
\end{align*}
We now note that $\dot{\phi}$ is a Jacobi field, and as in \cite{marx2021variations}, section \S 5, $\dot{\phi}$ admits an expansion (thought of as a function on $\Gamma = \gamma \times [0, \rho)$) which is odd up to order $2m^* - 1$:
\begin{align} \label{dotPhiExpansion}
\dot{\phi}(s,\tx) &= \frac{\dot{\phi}_0(s)}{\tx} + \tx \dot{\phi}_2(s) + (\text{odd}) + \tx^{2m^*-1} \dot{\phi}_{2m^*}(s) + \tx^{2m^*} \dot{\phi}_{2m^*+1}(s) + O(\tx^{2m^*+1}) 
\end{align}
Using that $d\tx(\nu) = \tx u_{\tx} + O(x^3)$ is even up to order $2m^*$ (see equation \eqref{NormalExpansion}), 
\begin{align*}
\dot{\tx}_{\dot{\phi}} &= \dot{\phi} d\tx(\nu) \\
&= O(\tx) \\
\tr_{Y_{\gamma}}(\omega_2) &= O(\tx^2) \\
\bJac_{Y_{\gamma}} &= O(1) 
\end{align*}
Furthermore, $\dot{\tx}_{\dot{\phi}}$ is odd up to order $2m^* - 1$, and $\tr_{Y_{\gamma}}(\omega_2)$ is even up to order $2m^* + 2$. In particular 
\[
\implies [ \dot{\tx}_{\phi} \tr_{Y_{\gamma}}(\omega_2) \bJac_{Y_{\gamma}}]^{(2m^*)} = 0
\]
\subsubsection{$B$} \label{B1Section}
We start by writing
\begin{align*}
P_t^*[\tr_{Y_{\gamma_t}}(\omega_2)] &= h(t)^{\alpha \beta} \omega_2(v_{\alpha}(t), v_{\beta}(t))
\end{align*}
here $\{v_{\alpha}(t) = v_{\alpha}(t, P_t(p))\}$ is a frame for $Y_{\gamma_t}$ achieved by pushing forward a frame on $Y_{\gamma}$ to $Y_{\gamma_t}$ by $(P_t)_*$. Similarly $h(t)_{\alpha \beta} = g(v_{\alpha}(t), v_{\beta}(t))$ is the induced metric on $Y_{\gamma_t}$. We differentiate this
\begin{align*}
\frac{d}{dt} P_t^*[\tr_{Y_{\gamma_t}}(\omega_2)] &= \frac{d}{dt} h(t)^{\alpha \beta} \omega_2(v_{\alpha}(t), v_{\beta}(t))	\\
&= \dot{h}^{\alpha \beta} \omega_2(v_{\alpha}, v_{\beta}) + h^{\alpha \beta} (\n_{\dot{\phi} \nu} \omega_2)(\va, \vb) + h^{\alpha\beta} [\omega_2( \dot{\va}, \vb) + \omega_2(\va, \dot{\vb})] \\
&= h^{\alpha \beta} (\n_{\dot{\phi} \nu} \omega_2)(\va, \vb) \\
&= \dot{\phi} h^{\alpha \beta} (\n_{\nu} \omega_2)(\va, \vb)
\end{align*}
We compute
\begin{align} \label{OmegaCovarDiff}
\nabla_{\nu}  \omega_2 &= \n_{\nu} \left( \frac{\omega_{2}(p)}{\tx^2} \tx^2 \right) \\ \nonumber
&= 2x \nu(x) \frac{\omega^2_{ij}(p) d \ty^i d \ty^j}{\tx^2} + \tx^2 \n_{\nu} \left(\frac{\omega^2_{ij}(p) d \ty^i d \ty^j}{\tx^2}\right) \\ \nonumber
&= \frac{2}{\tx} \nu(\tx) \omega_2 \\
&= 2 x \omega_2
\end{align}
having noted that $\n_{\nu} \left(\frac{\omega^2_{ij} d \ty^i d \ty^j}{\tx^2}\right) = 0$ since the basis forms $\frac{d \ty^i}{\tx}$ are parallel with respect to the hyperbolic metric. We also noted that $\nu(\tx) = \tx^2$ coming from parameterization in section \S \ref{HemisphereParamSection}. Thus
\begin{align*}
B& = \FPz \int_{Y_{\gamma}} \tx^z \left(\frac{d}{dt} P_t^*[\tr_{Y_{\gamma_t}}(\omega_2)] \right) dA_{Y_{\gamma}} \\
&= 2 \FPz \int_{Y_{\gamma}} \tx^{z+1} \dot{\phi} \tr_{Y_{\gamma}}(\omega_2) dA_{Y_{\gamma}} 
\end{align*}
\subsubsection{Recovering $\omega_2$}
Adding our terms together, we've shown that
\begin{align*}
\frac{d}{dt} \frac{d^2}{d \delta^2} \RA(Y_{\gamma_t(y), \delta}, g_{\delta}) \Big|_{\delta = 0} \Big|_{t = 0} &= 2 \FPz \int_{Y_{\gamma}} \tx^{z+1} \dot{\phi} \tr_{Y_{\gamma}}(\omega_2) dA_{Y_{\gamma}}
\end{align*}
To recover $\omega_2$, we explicitly compute $\dot{\phi}$ using that our underlying minimal surface is $HS^{2m^*}$. In the appendix section \S \ref{SolvingJacobiGeneral}, we show that we can find particular Jacobi fields using separation of variables
\[
\dot{\phi}(\theta_{2m^*-1}, \tx) = f_2(\theta_{2m^*-1}) \frac{\overline{B}_{2m^*, 2}(\tx)}{\tx}
\] 
Let $f_2$ be
\[
f_2(\theta_{2m^*-1}) = \cos^2 \phi - \frac{1}{2m^* - 1} \sin^2 \phi
\]
where we parameterize
\[
S^{2m^*-1} = \{(y_1, y_2, \dots, y_{2m}) = (\cos \phi, \sin \phi F_{2m^*-2}(\theta_{2m^*-2}))\}
\]
where $F_{2m^*-2}$ is a parameterization of $S^{2m^*-2}$. As shown in the appendix \S \ref{SolvingJacobiGeneral}, equation \eqref{B2Equation}:
\begin{align*}
\overline{B}_{2m^*, 2}(\tx) &= \frac{\F21(-m^*, -1/2 - m^*, 1/2-m^*; \tx^2) + c \tx^{2m^* + 1}}{(1-\tx^2)^{m^*}} \\
c & = - \F21(-m^*, -1/2 - m^*, 1/2 - m^*, 1) = - \frac{\Gamma(1+m^*) \Gamma(1/2-m^*)}{\Gamma(1/2)}
\end{align*}
where $\F21(a,b,c;x)$ is a hypergeometric function and our choice of $c$ means that $\overline{B}_{2m^*,2}$ is defined on all of $[0,1]$. We can compute
\begin{align*}
\frac{d}{dt} \frac{d^2}{d \delta^2} \RA(Y_{\gamma_t(y), \delta}, g_{\delta}) \Big|_{\delta = 0} \Big|_{t = 0}(\dot{\phi}) &= 2\FPz \int_{HS^m} \tx^{z+1} \dot{\phi} \tr_{HS^m}(\omega_2) dA_{HS^m} \\
&= 2\FPz \int_{x = 0}^1 \int_{S^{2m^*-1}} \tx^{z} f_2(\theta_{2m^*-1}) \overline{B}_{2m^*,2}(\tx) \left[ h^{\tx\tx} \omega_2(v_{\tx}, v_{\tx}) + h^{ij}\omega_2(v_i, v_j) \right]
\end{align*}
Again, using our parameterization from \S \ref{HemisphereParamSection}, we have 
\begin{align*}
h^{\tx\tx} \omega_2(v_{\tx}, v_{\tx}) &= \tx^4 \omega_2(F_{2m^*-1}(\theta_{2m^*-1}), F_{2m^*-1}(\theta_{2m^*-1})) \\
h^{ij} \omega_2(v_i, v_j) &= \tx^2 \tr_{S^{2m^*-1}}(\omega_2)(\theta_{2m^*-1}) \\
dA_{HS^{2m^*}} &= \frac{(1-\tx^2)^{m^*-1}}{\tx^{2m^*}} d\tx dA_{S^{2m^*-1}}(\theta_{2m^*-1})
\end{align*}
So that our main computation of interest becomes
\begin{align*}
\frac{d}{dt} \frac{d^2}{d \delta^2} \RA(Y_{\gamma_t(y), \delta}, g_{\delta}) \Big|_{\delta = 0} \Big|_{t = 0}(\dot{\phi}) &= 2\FPz \int_{x = 0}^1 \int_{S^{2m^*-1}} \tx^{z} f_2(\theta_{2m^*-1}) \overline{B}_{2m^*}(\tx) \Big[ \tx^4 \omega_2(F_{2m^*-1}(\theta_{2m^*-1}), F_{2m^*-1}(\theta_{2m^*-1})) \\
& \qquad \qquad \qquad \qquad + \tx^2 \tr_{S^{2m^*-1}}(\omega_2)(\theta_{2m^*-1})\Big] \frac{(1-\tx^2)^{2m^*-1}}{\tx^{2m^*}} \; dx dA_{S^{2m^*-1}}(\theta_{2m^*-1}) 
\end{align*}
In the appendix \S \ref{AngularComputation}, we compute
\begin{align*}
\int_{S^{2m^*-1}} f_2(\theta_{2m^*-1}) \omega_2(F_{2m^*-1}(\theta_{2m^*-1}), F_{2m^*-1}(\theta_{2m^*-1})) \; dA_{S^{2m^*-1}}(\theta_{2m^*-1}) &= \alpha \omega^2_{11} - \beta [ \omega^2_{22} + \omega^2_{33} + \dots + \omega^2_{(2m^*) (2m^*)}] \\
\int_{S^{2m^*-1}} f_2(\theta_{2m^*-1}) \tr_{S^{2m^*-1}}(\omega_2)(\theta_{2m^*-1}) \; dA_{S^{2m^*-1}}(\theta_{2m^*-1}) &= -\alpha \omega^2_{11} + \beta [ \omega^2_{22} + \omega^2_{33} + \dots + \omega^2_{(2m^*) (2m^*)}]
\end{align*}
where $\alpha, \beta > 0$ are explicitly computable in terms of $m^*$. As a result,
\begin{align} \nonumber
\frac{d}{dt} \frac{d^2}{d \delta^2} \RA(Y_{\gamma_t(y), \delta}, g_{\delta}) \Big|_{\delta = 0} \Big|_{t = 0}(\dot{\phi}) &= 2\FPz \int_0^1 \tx^{z+2} (1-x^2) \overline{B}_{2m^*,2}(\tx) [\alpha \omega^2_{11} - \beta (\omega^2_{22} + \dots + \omega^2_{(2m^*)(2m^*)})] \frac{(1-\tx^2)^{m^*-1}}{\tx^{2m^*}} dx \\ \label{DetectDifference}
&= 2 [\alpha \omega^2_{11} - \beta (\omega^2_{22} + \dots + \omega^2_{(2m^*)(2m^*)})] \cdot I_2 \\ \nonumber
I_2 &:= \FPz \int_0^1 \tx^{z-2m^*} \overline{B}_{2m^*,2}(x) (1-\tx^2)^{m^*} \\ \nonumber
&= \FPz \int_0^1 \tx^{z-2m^*} \left[ \F21(-m^*, -1/2 - m^*, 1/2-m^*; \tx^2) + c \tx^{2m^* + 1} \right]
\end{align}
In the appendix \S \ref{FinitePartComputation}, we compute that $I_2 \neq 0$, and in general
\begin{align*}
I_{k, m^*} &= \FPz \int_0^1 x^{z-2m^* - 2k} \left[ \F21(-m^*, -1/2 - m^*, 1/2-m^*; x^2) + c x^{2m^* + 1} \right] \neq 0
\end{align*}
From this, we see that
\begin{align*}
\frac{1}{2I_2 \alpha} \frac{d}{dt} \frac{d^2}{d \delta^2} \RA(Y_{\gamma_t(y), \delta}, g_{\delta}) \Big|_{\delta = 0} \Big|_{t = 0} &= \omega^2_{11} - (\beta/\alpha)[\omega^2_{22} + \dots + \omega^2_{(2m^*)(2m^*)}]
\end{align*}
Using this and lemma \ref{TraceRecoveryLemma}, we see that we can recover $\omega^2_{11}$, and by symmetry each $\omega^2_{kk}$. Changing coordinates
\begin{align*}
y_1 \to \ty_1 &= \frac{y_1 + y_2}{\sqrt{2}} \\
y_2 \to \ty_2 &= \frac{y_1 - y_2}{\sqrt{2}}	
\end{align*}
we see that $\omega_0(p)$ is still orthonormal with respect to $\{\ty_1, \ty_2, y_3, \dots, y_{2m^*}\}$. Repeating the proof with this new coordinate basis, we can recover
\[
\omega_{2}(\p_{y_1} + \p_{y_2}, \p_{y_1} + \p_{y_2})
\]
and hence recover $\omega_2(\p_{y_1}, \p_{y_2})$. Repeating this argument for every other pair of indices, we recover all of $\omega_2$. \qed
\subsection{Codimension $1$, Induction $k \leq 2m^*-2$}
We assume knowledge of $\omega_0, \dots, \omega_{k-1}$ (for $k \leq 2m^*-2$, even) and then differentiate $k$ times to isolate the $\omega_k$ term. As in the base case, we vary the boundary curve to induce an anisotropy which allows us to isolate $\omega_{k,ii}$ for each $i$. \nl \nl
We first apply $\frac{d^k}{d \delta^k}$. Arguing as in the base case:
\begin{align*}
\frac{d^k}{d \delta^k} \RV(Y_{\gamma, \delta}, g_{\delta}) &= \FPz \int_{\tY_{\gamma, \delta}} \frac{d^k}{d \delta^k} \left[ P_{\delta}^*(\tx)^z P_{\delta}^*(dA_{Y_{\gamma(\ty), \delta}}) \right] \\ 
&= \FPz \Big[\int_{Y_{\gamma}} z x^{z-1} \p_{\delta}^k(P_{\delta}^*(\tx)) dA_{Y_{\gamma}} \\
& \qquad + \int_{Y_{\gamma}} x^z \left(\frac{d^k}{d \delta^k} P_{\delta}^*(dA_{\tY_{\gamma, \delta}}) + \sum_{\ell = 1}^{k-1} c_{\ell, k} \left[\frac{d^{\ell}}{d \delta^{\ell}} P_{\delta}^*(\tx^z) \right]  \left[\frac{d^{k-\ell}}{d \delta^{k-\ell}} P_{\delta}^*\left(dA_{\tY_{\gamma, \delta}} \right) \right] \right) \Big] \\
&= \FPz \left[\int_{Y_{\gamma}} z \tx^{z-1} \p_{\delta}^k(P_{\delta}^*(\tx)) dA_{Y_{\gamma}} + \int_{Y_{\gamma}} \tx^z \frac{d^k}{d \delta^k} P_{\delta}^*(dA_{\tY_{\gamma, \delta}}) \right] + K(\omega_0, \omega_2, \dots, \omega_{k-1})  
\end{align*}
Here, we've noted that all of the $\delta$ derivatives which act on both $P_{\delta}^*(\tx^z)$ and $P_{\delta}^*(dA_{\tY_{\gamma, \delta}})$ will depend on at most $k-1$ derivatives of $g_{\delta}$, and hence be determined since we assume knowledge of $\omega_0, \dots, \omega_{k-1}$ (which determines $\p_{\delta}^{\ell} g_{\delta}$ for any $\ell \leq k - 1$).  Mimicking the computation for $k = 2$, we have
\begin{align*}
\FPz \int_{Y_{\gamma}} z \tx^{z-1} \p_{\delta}^k(P_{\delta}^*(\tx)) dA_{Y_{\gamma}} &= \int_{\gamma} [d\tx(\p_{\delta}^k P_{\delta}) \sqrt{\det \bh}]^{(2m^*)} dA_{\gamma}
\end{align*}
Similarly, as in the $k = 2$ case subsection \S \ref{DTermAnalysis}, one can show that 
\begin{align*}
\frac{d^k}{d \delta^k} P_{\delta}^*(dA_{\tY_{\gamma, \delta}}) \Big|_{\delta = 0} &= \div_{\tY_{\gamma}}(\p_{\delta}^kP_{\delta}) + \tr_{\tY_{\gamma}}(\dot{g}^{(k)}) + K(\omega_0, \omega_2, \dots, \omega_{k-1}) \\
&= \div_{\tY_{\gamma}}(\p_{\delta}^kP_{\delta}) + (k-2)!\tr_{\tY_{\gamma}}(x^{k-2} \omega_k) + K(\omega_0, \omega_2, \dots, \omega_{k-1})
\end{align*}
Having recalled equation \eqref{NthMetricDerivNoLog},
\begin{align} 
k \leq 2m^* - 2 \implies \dot{g}^{(k)} &= \partial_{\delta}^{(k)} g_{\delta} \\
&= (k-2)! \tx^{k-2} \omega_k + O(\omega_0, \dots, \omega_{k-1})
\end{align}
%
This yields:
\begin{align} \label{HigherOrderDeltaDeriv}
\frac{d^k}{d \delta^k} \RA(Y_{\gamma_t(y), \delta}, g_{\delta}) \Big|_{\delta = 0} &= (k-2)!\FPz \int_{Y_{\gamma_t}} \tx^z \tr_{Y_{\gamma_t}}(\tx^{k-2} \omega_k) dA_{Y_{\gamma_t}} + K(\omega_0, \dots, \omega_{k-1}) 
\end{align}
As in the base case, we now show that we can recover the trace of $\omega_k$ from the above 
\begin{lemma} \label{TraceRecoveryLemmaGeneralk}
	Equation \eqref{HigherOrderDeltaDeriv} recovers the trace of $\omega_k$ at $p$ for any $k \geq 2$:	
\end{lemma}
\noindent \Pf We set $\gamma(\ty) = S^{2m^*-1}$ and $Y = HS^{2m^*}$ in equation \eqref{HigherOrderDeltaDeriv}. Mimicking the computations exactly as in lemma \ref{TraceRecoveryLemma}, we have
\begin{align} \label{TracePrelimEqnGeneralK}
	\frac{1}{(k-2)!} \frac{d^k}{d \delta^k} \RA(Y_{S^{2m^*-1}, \delta}, g_{\delta}) \Big|_{\delta = 0} &= \FPz \int_{HS^{2m^*}} \tx^{z + k - 2} \tr_{HS^{2m^*}}(\omega_k) dA_{HS^{2m^*}} \\ \nonumber
	&= [\omega^k_{11} + \dots + \omega^k_{(2m^*, 2m^*)}] \left( c_1 I_1 + c_2 I_2 \right)
\end{align}
where 
\begin{align*}
	c_1 [\omega^k_{11} + \dots + \omega^k_{(2m^*, 2m^*)}] &= \int_{S^{2m^*-1}}\omega_k(F_{2m^*-1}(\theta), F_{2m^*-1}(\theta)) dA_{S^{2m^*-1}} \\
	c_2 [\omega^k_{11} + \dots + \omega^k_{(2m^*, 2m^*)}] &= \int_{S^{2m^*-1}} \omega_k(F_{2m^*-1}(\theta), F_{2m^*-1}(\theta)) dA_{S^{2m^*-1}} \\
	I_1 &= \FPz \int_0^1 \tx^{z+k+4-2m^*}(1-x^2)^{m^*-1} \\
	I_2 &= \FPz \int_0^1 \tx^{z+k+2-2m^*}(1-x^2)^{m^*-1}
\end{align*}
Exactly as in lemma \ref{TraceRecoveryLemma}, we have
\begin{align*}
c_1 &= \frac{1}{2m^*-1} \left[\int_{\phi = 0}^{\pi} (\sin \phi)^{2m^*}d\phi\right] \text{Vol}(S^{2m^*-2}) \\
c_2 &= (2m^*-1) c_1
\end{align*}
We also mimic the computation for $I_1, I_2$ via integration by parts:
\begin{align*}
	I_1 &= \FPz \int_0^1 x^{z+k+4-2m^*}(1-x^2)^{m^*-1} \\
	&= \FPz \frac{2(m^*-1)}{z+k+5-2m^*} \int_0^1 x^{z+k+6-2m^*}(1-x^2)^{m^*-2} \\
	&= \FPz \frac{2^{m^*-1}(m^*-1)!}{(z+k+5-2m^*)(z+k+7-2m^*)\cdots(z+k+1)} \int x^{z+k+2}  \\
	&= \frac{2^{m^*-1} (m^*-1)!}{(k+5-2m^*)(k+7-2m^*) \cdots (k -1) (k+ 1)} \cdot \frac{1}{k+3}
\end{align*}
Similarly
\begin{align*}
	I_2 &= \FPz \int_0^1 x^{z+k+2-2m^*}(1-x^2)^{m^*-1} \\
	&= \FPz \frac{2(m^*-1)}{(z+k+3-2m^*)} \int_0^1 x^{z+4-2m^*}(1-x^2)^{m^*-2} \\
	&= \FPz \frac{2^{m^*-1} (m^* - 1)!}{(z + k + 3 - 2m^*) (z + k + 5 - 2m^*) \cdots (z + k -1)} \int_0^1 x^{z + k} \\
	&= \frac{2^{m^* - 1} (m^* - 1)!}{(k + 3 - 2m^*)(k + 5 - 2m^*) \cdots (k-1) (k + 1)} \\
	&= \frac{k + 3}{2m^* + k -3} I_1
\end{align*}
Together, this gives 
\begin{align*}
	c_1 I_1 &= \frac{c_2}{2m^*-1} \cdot \left[ (-1) \frac{3-2m^* - k}{k+ 3} \right] I_2 \\
	&= c_2 I_2 \frac{2m^* +k - 3}{(2m^*-1) (k+3)} \\
	\implies c_1 I_1 + c_2 I_2 &= c_2 I_2 \left[1 + \frac{2m^* +k - 3}{(2m^*-1) (k+3)} \right] \\
	&\neq 0
\end{align*}
This computation along with equation \eqref{TracePrelimEqnGeneralK} finishes the proof. Note that if $k \geq 2m^*$, there is no need to compute $\FPz$ as the integrals
\begin{align*}
f_1(z) &= \int_0^1 x^{z+k+4-2m^*}(1-x^2)^{m^*-1} \\
f_2(z) &= \int_0^1 x^{z+k+2-2m^*}(1-x^2)^{m^*-1}
\end{align*}
are both holomorphic near $z = 0$ when $k \geq 2m^*$. \qed \nl \nl
We now vary the $k$th $\delta$-derivatives of renormalized volume along a family of curves and define
\begin{align*}
L(t) &:= \FPz \int_{Y_{\gamma_t}} \tr_{Y_{\gamma_t}}(\tx^{z+ k-2} \omega_k) dA_{Y_{\gamma_t}}
\end{align*}
Equation \eqref{HigherOrderDeltaDeriv} tells us that $\frac{d^k}{d \delta^k} \RA(Y_{\gamma_t(y), \delta}, g_{\delta}) \Big|_{\delta = 0}$ and knowledge of $\{\omega_0, \dots, \omega_{k-1}\}$ immediately determines $L(t)$, so it suffices to compute variations of $L(t)$. As in the base case, we now consider $\gamma_t$ to be a perturbation over $\gamma_0 = S^{2m^*-1}$ with $Y_0 = HS^{2m^*}$ and write
\[
Y_t = P_t(HS^{2m^*}) = \{\exp_{HS^{2m^*}}(\phi_t(p) \nu(p)) \}
\]
to compute
\begin{align*}
L(t) &= \FPz \int_{HS^{2m^*}} P_t^*(\tx)^{z + k-2} P_t^*(\tr_{Y_t} \omega_k) P_t^*(dA_{Y_t}) \\
\frac{d}{dt} L(t) \Big|_{t = 0} &= \FPz \int_{HS^{2m^*}} (z + k-2) \tx^{z + k-3} \dot{\tx}_{\phi} \tr_{HS^{2m^*}} \omega_k dA_{HS^{2m^*}} + \tx^{z + k-2} \frac{d}{dt} P_t^*(\tr_{Y_t} \omega_k) dA_{HS^{2m^*}} \\
&= A + B
\end{align*}
As before, we compute in pieces
\subsubsection{$A$}
We have for $k \geq 3$
\begin{align*}
A&= \FPz \int_{HS^{2m^*}} (z + k-2) \dot{\tx}_{\phi} \tx^{z + k-3} \tr_{HS^{2m^*}} \omega_k dA_{HS^{2m^*}} \\
&= \FPz z \int_{HS^{2m^*}} \dot{\phi} \nu(\tx) \tx^{z + k-3} \tr_{HS^{2m^*}} \omega_k dA_{HS^{2m^*}} \\
& \quad +  (k-2) \FPz \int_{HS^{2m^*}}  \dot{\phi} \nu(\tx) \tx^{z + k-3} \tr_{HS^{2m^*}}(\omega_k) dA_{HS^{2m^*}} \\
&= \int_{S^{2m*-1}} [\tx^{k} \dot{\phi} (\tr_{HS^{2m^*}} \omega_k) \sqrt{\det \bh_{HS^{2m^*}}}]^{2m^*}   \\
& \quad + (k-2) \FPz \int_{HS^{2m^*}}  \dot{\phi} \tx^{z + k-1} \tr_{HS^{2m^*}}(\omega_k) dA_{HS^{2m^*}} 
\end{align*}
here, we've localized one of the terms using techniques from \S \ref{mechanics} and used that $\nu(x) = x^2$ for $Y = HS^{2m^*}$. We now note that $k$ is even, and so
\[
\tx^k (\tr_{HS^{2m^*}} \omega_k) \sqrt{\det \bh_{HS^{2m^*}}} = O(x^{2 + k})
\]
and is even to order $2m^* + 2 + k$. Furthermore, $\dot{\phi}$ is $O(x^{-1})$ and odd to order $2m^* - 1$ (see equation \eqref{dotPhiExpansion}). As a result, the integrand, $\tx^{k} \dot{\phi} (\tr_{HS^{2m^*}} \omega_k) \sqrt{\det \bh_{HS^{2m^*}}}$ is $O(x^{1+k})$ and odd to order $2m^* + 1 + k$. In particular, there is no $x^{2m^*}$ coefficient so 
\begin{align*}
[\tx^{k} \dot{\phi} (\tr_{HS^{2m^*}} \omega_k) \sqrt{\det \bh_{HS^{2m^*}}}]^{2m^*} &= 0 \\
\implies A &= (k-2) \FPz \int_{HS^{2m^*}}  \dot{\phi} \tx^{z + k-1} \tr_{HS^{2m^*}}(\omega_k) dA_{HS^{2m^*}} 
\end{align*}
We also remark that 
\[
k \geq 2m^* \implies [\tx^{k} \dot{\phi} (\tr_{HS^{2m^*}} \omega_k) \sqrt{\det \bh_{HS^{2m^*}}}]^{2m^*} = 0
\]
independent of the parity of $k$.
\subsubsection{$B$}
As in the base case with section \S \ref{B1Section}, we have
\begin{align*}
\frac{d}{dt} P_t^*(\tr_{Y_t} \omega_k)&= h^{\alpha \beta} \dot{\phi} (\n_{\nu} \omega_k)(\va, \vb) \\
&= 2\tx \dot{\phi} \tr_{HS^{2m^*}}(\omega_k) \\
\implies B &= 2  \FPz \int_{HS^{2m^*}} \dot{\phi} \tx^{z + k-1} \tr_{HS^{2m^*}}(\omega_k) dA_{HS^{2m^*}}
\end{align*}
%
%
\subsubsection{In Sum, induction} \label{GeneralMProof}
Together, we have
\begin{align*}
\frac{d}{dt} L(t)\Big|_{t = 0} &= A + B\\
& = k \; \FPz \int_{HS^{2m^*}} \dot{\phi} \tx^{z+k-1} \tr_{HS^{2m^*}}(\omega_k) dA_{HS^{2m^*}}
\end{align*}
Recall from \S \ref{HemisphereParamSection}
\begin{align*}
dA_{HS^{2m^*}} &= \frac{(1-\tx^2)^{m^*-1}}{\tx^{2m^*}} d\tx dA_{S^{2m^*-1}} 
\end{align*}
We can then compute
\begin{align} \label{DifferentiateEqn}
\frac{d}{dt} L(t)\Big|_{t = 0} &= k \FPz \int_{HS^{2m^*}} \dot{\phi} \tx^{z+k-1 - 2m^*}(1-\tx^2)^{m^*-1} \tr_{HS^{2m^*}}(\omega_k) dA_{S^{2m^*-1}} d\tx \\ \nonumber
\tr_{HS^{2m^*}}(\omega_k) &= \tx^2 \tr_{S^{2m^*-1}}(\omega_k)(\theta) + \tx^4 \omega_k(F_{2m^*-1}(\theta), F_{2m^*-1}(\theta)) \\ \nonumber
dA_{S^{2m^*-1}} &= (\sin^{2m^*-2} \phi) d \phi dA_{S^{2m^*-2}}(\theta)
\end{align}
Here, $F_{2m^*-1}(\theta)$ is the position vector on $S^{2m^*-1}$ and we've again parameterized $S^{2m^*-1}$ in polar coordinates. Using the parameterization according to section \S \ref{HemisphereParamSection}, we can write
\begin{align*}
	\tr_{HS^m}(\omega_k) &= h^{\tx\tx} \omega_k(v_{\tx}, v_{\tx}) + h^{ij} \omega_k(v_i, v_j) \\
	&= \tx^4 \omega_k(F_{2m^*-1}(\theta), F_{2m^*-1}(\theta)) + \tx^2 \bh^{ij}(\theta)\omega_k(F_{2m^*-1,i}(\theta), F_{2m^*-1,j}(\theta)) \\
	&= \tx^4 \omega_k(F_{2m^*-1}(\theta), F_{2m^*-1}(\theta)) + \tx^2 \tr_{S^{2m^*-1}}(\omega_k)
\end{align*}
According to section \S \ref{SolvingJacobiGeneral}, we can again choose the following Jacobi field 
%
\begin{align*}
	\dot{\psi} &= \frac{\overline{B}_{2,2m^*}(\tx)}{\tx(1-\tx^2)^{m^*}} \left[\cos^2 \phi - \frac{1}{m-1} \sin^2 \phi \right] \\
	&= \frac{\overline{B}_{2,2m^*}(\tx)}{\tx(1-\tx^2)^{m^*}} \cdot f_2(\phi) 
\end{align*}
So that our integral becomes
\begin{align*}
\frac{d}{dt} L(t)\Big|_{t = 0} &= k \FPz \int_{HS^{2m^*}} \dot{\psi} \tx^{z+k-1 - 2m^*}(1-\tx^2)^{m^*-1} \tr_{HS^{2m^*}}(\omega_k) dA_{S^{2m^*-1}} d\tx \\
&= k \FPz \int_0^1 \tx^{z+k-2m^*} \overline{B}_{2,2m^*}(x) (1-\tx^2)^{-1} d\tx \int_{S^{2m^*-1}}f_2(\theta)\tr_{S^{2m^*-1}}(\omega_k)(\theta) dA_{S^{2m^*-1}}(\theta) \\
& + k \FPz \int_0^1 \tx^{z+k-2m^*+2} \bB_{2,2m^*}(\tx)  (1-\tx^2)^{-1} dx \int_{S^{2m^*-1}} f_2(\theta)\omega_k(F_{2m^*-1}(\theta), F_{2m^*-1}(\theta)) dA_{S^{2m^*-1}}(\theta) 
\end{align*}
As in the base case, we show compute the angular integrals in the appendix \S \ref{AngularComputation}
\begin{align} \label{I1Integral}
I_1 &= \int_{S^{2m^*-1}} f_2(\theta) \omega_k(F_{2m^*-1}(\theta_{2m^*-1}), F_{2m^*-1}(\theta_{2m^*-1})) dA_{S^{2m^*-1}}(\theta_{2m^*-1}) \\ \nonumber
&= \left(\alpha\omega^k_{11} - \beta[\omega^k_{22} + \dots + \omega^k_{(2m^*)(2m^*)}] \right)\left[ (2/(2m^*-1)) \int_0^{\pi} \cos^2(\phi) \sin^{2m^*}(\phi) \right] \\ \nonumber
&= \left(\alpha\omega^k_{11} - \beta[\omega^k_{22} + \dots + \omega^k_{(2m^*)(2m^*)}] \right)c_{m^*} \\ \label{I2Integral}
I_2 &= \int_{S^{2m^*-1}} f_2(\theta)\omega_k(F_{2m^*-1}(\theta), F_{2m^*-1}(\theta)) dA_{S^{2m^*-1}}(\theta)  \\ \nonumber
&= - I_1
\end{align}
where $c_{m^*}$ is non-zero and $\alpha, \beta > 0$ for all $m^*$. Our overall integral becomes
\begin{align*}
\frac{d}{dt} L(t)\Big|_{t = 0} &=  k c_{m^*} \left(\alpha\omega^k_{11} - \beta[\omega^k_{22} + \dots + \omega^k_{(2m^*)(2m^*)}] \right) \FPz \int_0^1 x^{z+k-2m^*} \bB_{2,2m^*}(x) dx \\
&=  k c_{m^*} \left(\alpha\omega^k_{11} - \beta[\omega^k_{22} + \dots + \omega^k_{(2m^*)(2m^*)}] \right) I_{k,m^*} \\
I_{k,m^*} &:= \FPz \int_0^1 x^{z+k-2m^*} \bB_{2,2m^*}(x) dx
\end{align*}
In the appendix \S \ref{FinitePartComputation}, we show that $I_{k,m^*}$ is non-zero for all $k, m^*$, so that 
\[
\frac{d}{dt} L(t)\Big|_{t = 0} =  C(k, m^*) \left(\alpha\omega^k_{11} - \beta[\omega^k_{22} + \dots + \omega^k_{(2m^*)(2m^*)}] \right)
\]
where $C(k, m^*) \neq 0$. This computation, in tandem with the trace recovery result from lemma \ref{TraceRecoveryLemmaGeneralk} allows us to determine $\omega_{k,11}$. Repeating this for any orthonormal basis for $\omega_0(p)$ determines $\omega_k(p)$ in full. \qed
\subsubsection{$k \geq 2m^*$}
When $k \geq 2m^*$, the asymptotic expansion of $g$ contains factors of $x^{k} \log(x)^p$, since we've assumed that $g$ is partially even to only order $2m^*$. To resolve this, we first determine $\omega^{k P_k}$ where we recall that $x^k\log(x)^{P_k}$ is the highest order $\log(x)$ term associated to the power $x^k$. We then determine the lower order coefficients $\omega^{k (P_k - \ell)}$ for $1 \leq \ell \leq P_k$ inductively. We then induct again and determine $\omega^{k b}$ for higher orders of $k$. \nl 
\indent To begin, we recall equation \eqref{FPDelta}
\begin{align*}
\FPd \frac{1}{\log(\delta)^a} \frac{d^b}{d \delta^b} \delta^{c} \log(\delta)^d \Big|_{\delta = 0} &= \begin{cases}
	b! & d = a, \; c = b \\
	0 & d < a, \; c = b \\
	0 & c \neq b 
\end{cases} \\ \nonumber
\implies \FPd \frac{1}{\log(\delta)^a} \frac{d^b}{d \delta^b} g_{\delta} &= b! \tx^{b-2} \omega_{ba} + F(\{\tx^{b-2}\omega^{b \ell}_{\ell > a}\})
\end{align*}
With this, we can mimic much of the induction as to when there are no $\log(\delta)$ terms. We first compute
\begin{align*}
\FPd\frac{1}{\log(\delta)^{P_k}}\frac{d^k}{d \delta^k} \RV(Y_{\gamma, \delta}, g_{\delta}) &= \FPz \int_{\tY_{\gamma, \delta}} \FPd\frac{1}{\log(\delta)^{P_k}}\frac{d^k}{d \delta^k} \left[ P_{\delta}^*(\tx)^z P_{\delta}^*(dA_{Y_{\gamma(\ty), \delta}}) \right] \\ 
&= \FPz \Big[\int_{Y_{\gamma}} z \tx^{z-1} \left[\FPd\frac{1}{\log(\delta)^{P_k}}\p_{\delta}^k(P_{\delta}^*(\tx))\right] dA_{Y_{\gamma}} \\
& \qquad + \int_{Y_{\gamma}} \tx^z \FPd \frac{1}{\log(\delta)^{P_k}}\left(\frac{d^k}{d \delta^k} P_{\delta}^*(dA_{\tY_{\gamma, \delta}}) + \sum_{\ell = 1}^{k-1} c_{\ell, k} \left[\frac{d^{\ell}}{d \delta^{\ell}} P_{\delta}^*(\tx^z) \right]  \left[\frac{d^{k-\ell}}{d \delta^{k-\ell}} P_{\delta}^*\left(dA_{\tY_{\gamma, \delta}} \right) \right] \right) \Big] \\
&= \FPz \Big[\int_{Y_{\gamma}} z \tx^{z-1} \FPd \left( \frac{1}{\log(\delta)^{P_k}}\p_{\delta}^k(P_{\delta}^*(\tx)) \right) dA_{Y_{\gamma}} \\
& + \int_{Y_{\gamma}} \tx^z \FPd \left( \frac{1}{\log(\delta)^{P_k}} \frac{d^k}{d \delta^k} P_{\delta}^*(dA_{\tY_{\gamma, \delta}}) \right) \Big] \\
& + K(\omega_0, \omega_2, \dots, \{\omega_{k-1}^{\ell}\})  
\end{align*}
Now as in the $k < 2m^*$ case, we have
\begin{align} \label{Localization1}
\FPz \int_{Y_{\gamma}} z \tx^{z-1} \left( \FPd \frac{1}{\log(\delta)^{P_k}}\p_{\delta}^k(P_{\delta}^*(\tx)) \right) dA_{Y_{\gamma}} &= \int_{\gamma} \left[ \FPd \frac{1}{\log(\delta)^{P_k}}\p_{\delta}^k(P_{\delta}^*(\tx)) \sqrt{\det \bh}\right]^{(2m^*)} \\ \nonumber
&= \int_{\gamma} \left[ \FPd \frac{1}{\log(\delta)^{P_k}} d\tx\left(\p_{\delta}^k(P_{\delta})\right) \sqrt{\det \bh}\right]^{(2m^*)}
\end{align}
And similarly
\begin{align*}
\FPd \frac{1}{\log(\delta)^{P_k}} \frac{d^k}{d \delta^k} P_{\delta}^*(dA_{\tY_{\gamma, \delta}}) &= \div_{\tY_{\gamma}}\left( \FPd \frac{1}{\log(\delta)^{P_k}} \p_{\delta}^k P_{\delta} \right) + \tr_{\tY_{\gamma}} \left( \FPd \frac{1}{\log(\delta)^{P_k}} \p_{\delta}^k g_{\delta}\right)
\end{align*}
And as before
\begin{align*}
\FPz \int_{\tY_{\gamma}} \tx^z \div_{\tY_{\gamma}}\left( \FPd \frac{1}{\log(\delta)^{P_k}} \p_{\delta}^k P_{\delta} \right) &= - \FPz z \int_{\tY_{\gamma}} d\tx\left(\FPd \frac{1}{\log(\delta)^{P_k}} \p_{\delta}^k P_{\delta} \right)
\end{align*}
which cancels with the term from equation \eqref{Localization1}. We now use equation \eqref{FPDelta}:
\begin{align*}
\FPd \frac{1}{\log(\delta)^{P_k} } \p_{\delta}^k g_{\delta}  &= k! \tx^{k-2}\omega_{k P_k}
\end{align*}
noting that $P_k$ is the highest $\log(x)$ power associated to $\tx^k$. This gives
\begin{align*}
\FPd\frac{1}{\log(\delta)^{P_k}}\frac{d^k}{d \delta^k} \RV(Y_{\gamma, \delta}, g_{\delta})&= k! \FPz  \int_{\tY_{\gamma}} \tx^{z+k-2} \tr_{Y_{\gamma}}(\omega_{k P_k}) dA_{\tY_{\gamma}} + K(\omega_0, \omega_2, \dots, \{\omega_{(k-1)\ell}\}) 
\end{align*}
Now using lemma \ref{TraceRecoveryLemmaGeneralk}, the assumed knowledge of $\{\omega_0, \omega_2, \dots, \{\omega_{(k-1) \ell}\}\}$, and varying $\gamma \to \gamma_t$, we can determine $\omega^{k P_k}$ in the same way as with $k < 2m^*$. \nl \nl
We now induct down from $P_k \to P_k - \ell$. Assume knowledge of $\{\omega_0, \omega_2, \dots, \{\omega_{(k-1)\ell}\}, \{\omega_{k(P_k - b)}\}_{0 \leq b \leq \ell - 1}\}$. Then we now compute
\begin{align*}
\FPd\frac{1}{\log(\delta)^{P_k - \ell}}\frac{d^k}{d \delta^k} \RV(Y_{\gamma, \delta}, g_{\delta}) &= k! \FPz  \int_{Y_{\gamma}} \tx^z \tr_{\tY_{\gamma}}(\omega_{k (P_k - \ell)}) dA_{\tY_{\gamma}} + K(\omega_0, \omega_2, \dots, \{\omega_{k-1}^{\ell}\}, \{\omega_{k(P_k - b)}\}_{0 \leq b \leq \ell - 1}\}) 
\end{align*}
Using the same arguments as when determining $\omega^{k P_k}$ but noting that
\[
\FPd \frac{1}{\log(\delta)^a} \frac{d^b}{d \delta^b} g_{\delta} = b! \tx^{b-2} \omega_{ba} + F(\{\tx^{b-2}\omega_{b \ell}\}_{\ell > a})
\]
we can again use the same proof as in lemma \ref{TraceRecoveryLemmaGeneralk} to determine the trace and then varying $\gamma \to \gamma_t$, we can recover $\omega_{k \ell}$. Repeating this for all $\ell \geq 0$ inductively finishes the recovery at order $k$. \nl 
\indent Now repeating this for arbitrary $k > 2m^*$, we can determine all coefficients, $\{\omega_k, \omega_{ab}\}$ in the expansion given by equation \eqref{EvenExpansion}.

\subsection{Codimension $> 1$}
We consider $\RV(Y^{2m^*})$ for minimal surfaces in $(M^n, g)$ and $n \geq 2m^* + 2$ with $g$ being partially even of order $2m^*$. As before, we can coordinatize a neighborhood about $p \in \partial M$ as $\{x, y_1, \dots, y_{n-1}\}$. Choose any $2m^*$ coordinates among $\{y_1, \dots, y_{n-1}\}$ and WLOG label them $\{y_1, \dots, y_{2m^*}\}$. As before, let $\{\ty_1, \dots, \ty_{n-1}\}$ denote the coordinates for the blow up of $\overline{M}$ about $p$, namely $\H^{n}$. \nl \nl
We now choose the limiting curve, $\gamma^{2m^*-1} \subset \R^{2m^*}_{\ty_1, \dots, \ty_{2m^*}} \subset \partial \H^n$. Because of this, $\tY_{\gamma} \subseteq \H^{2m^* + 1} \subset \H^n$. The reader can then check that the codimension $1$ proof translates immediately considering variations restricted so that $\gamma_t \subseteq \R^{2m^*} \subset \partial \H^n$ with $\gamma_{t = 0} = \gamma$ given by:
\[
\gamma = \{\ty_1^2 + \dots + \ty_{2m^*}^2 = 1\}
\]
this will determine each $\{\omega^{k}(p), \omega^{ab}(p)\}$ restricted to the tangent plane, $T(p) = \text{span}\{\p_{y_1}, \dots, \p_{y_{2m^*}}\}\Big|_p$. Repeating this for every choice of $2m^*$ coordinates among $\{y_1, \dots, y_{n-1}\}$ determines $\Big\{\omega_k \Big|_{T_p M}, \omega_{ab} \Big|_{T_p M}\Big\}$. \nl \nl
\rmk \; When the codimension is $> 1$, we can actually prove the result without varying $\gamma \to \gamma_t$. Instead, lemma \ref{TraceRecoveryLemmaGeneralk} tells us that given $\omega_0$ globally on $\partial M$, we can determine $\tr_{T^{2m^*}}(\omega_2)$ for any $2m^*$-dimensional subspace $T^{2m^*} \subseteq T_p \partial M$. Since $2m^* < n-1$, this actually determines $\omega_2$ on the entirety of the tangent space by creating an invertible system of 
\[
\begin{pmatrix}
n - 1 \\
2m^*
\end{pmatrix} = \frac{(n-1)!}{(n-1-2m^*)! (2m^*)!} \geq n - 1
\]
equations with $n-1$ unknowns. If we assume that $\omega_2$ is diagonal (which we can by a change of basis since each $\{\omega_k, \omega_{ab}\}$ are symmetric), each equation will be the sum of $\omega^2_{ii}$ for $i$ in any subset of $\{1, \dots, n-1\}$ of size $2m^*$. Having determined $\omega_2 \Big|_{T_p M}$ for all $p$, we induct in the same way and determine all higher order tensor coefficients in the expansion of equation \eqref{EvenExpansion}. \qed 
\section{Applications} \label{Applications}
In this section we note some applications of theorem \ref{ExpansionThm}. Our first set of applications are corollaries \ref{SurfaceGGSURigidity} and \ref{GeneralGGSURigidity}, the corresponding rigidity results for CC partially even metrics using renormalized volume (compare to \cite{graham2019x}, Theorems 2 and 3)
\begin{corollary*} 
Suppose $(M, g)$, $(M, g')$ are two CC partially even metrics of order at least $2$. For $\gamma \subseteq \partial M$, let $Y^2_{\gamma, g}$ and $Y^2_{\gamma, g'}$ be the corresponding minimal surfaces (with respect to $g$ and $g'$ respectively) with $\partial Y^2_{\gamma, g} = \partial Y_{\gamma, g'}^2 = \gamma$. Suppose that $\RV_{g}(Y_{\gamma, g}^2) = \RV_{g'}(Y_{\gamma, g'}^2)$ for all $\gamma$. Then there exists a diffeomorphism $\psi: \overline{M} \to \overline{M}$ such that
\begin{align*}
\psi \Big|_{\partial M} &= Id \\
\psi^*(g') - g &= O(x^{\infty}) 
\end{align*}
Furthermore, if $g$ and $g'$ have log-analytic expansions and $\pi_1(\overline{M}, \partial \overline{M}) = 0$ then $\psi^*(g') = g$.
\end{corollary*}
\begin{corollary*} 
Suppose $(M, g)$, $(M, g')$ are two CC metrics that are partially even to order at least $m = 2m^*$ and $c(g) = c(g')$. For $\gamma^{m-1} \subseteq \partial M$, let $Y^m_{\gamma, g}$ and $Y^m_{\gamma, g'}$ be the corresponding minimal surfaces with $\partial Y^m_{\gamma, g} = \partial Y^m_{\gamma, g'} = \gamma$. Suppose that $\RV_{g}(Y_{\gamma, g}^m) = \RV_{g'}(Y_{\gamma, g'}^m)$ agree for all $\gamma$. Then there exists $\psi: \overline{M} \to \overline{M}$ such that
\begin{align*}
	\psi \Big|_{\partial M} &= Id \\
	\psi^*(g') - g &= O(x^{\infty}) 
\end{align*}
Furthermore, if $g$ and $g'$ have log-analytic expansions and $\pi_1(\overline{M}, \partial \overline{M}) = 0$ then $\psi^*(g') = g$. 
\end{corollary*}
\noindent \Pf We prove the $m = 2$ case first. Applying theorem \ref{ConformalInfinityThm}, the renormalized area functional on minimal surfaces determines the conformal infinity so $c(g) = c(g')$. Choose some $\omega_0 \in c(g)$. Using theorem \ref{ExpansionThm}, there exists chart maps $\psi_1$ and $\psi_2$ such that the metrics have log-analytic expansions, i.e.
\begin{align*}
\psi_1^*(g) &= \frac{dx^2 + \omega_0(y) + x^2 \omega_2(y) + \dots}{x^2} \\
\psi_2^*(g') &= \frac{d \rho^2 + \omega_0(y) + \rho^2 \omega_2(y) + \dots}{\rho^2}
\end{align*}
holds on some open neighborhood $U \supseteq \partial M$ for which $\psi_1, \psi_2$ are both defined. Here, $\rho$ and $x$ are special bdfs realizing the same $\omega_0$. Mapping $\rho \to x$ via $\psi_1 \circ \psi_2^{-1}$, we have
\begin{align*}
\implies g - (\psi_2 \circ \psi_1^{-1})^*(g') &= O(x^{\infty})
\end{align*}
Note that this analysis includes the presence of $x^k \log(x)^p$ terms.  \nl
\indent For full rigidity when $g$ and $g'$ are log-analytic, we follow the proof of theorem $3$ in \cite{graham2019x}. Let $\chi = \psi_2 \circ \psi_1^{-1}$. Log-analytic means that $x$ is real analytic on $\overline{M}$ and $g,\chi^*(g')$ have log-analytic (see \S \ref{polyHomAnalytic}) expansions in $x$ in a neighborhood of the boundary. This means $g = (\chi)^*(g')$ exactly on an open neighborhood $U \supseteq \partial M$. In particular, $U \supseteq [0, \eta)_x \times \partial M$ for some $\eta > 0$ by compactness of $\partial M$. \nl 
\indent As noted in \cite{graham2019x}, completeness of $M$ means that one can extend our isometry, $\chi$, via analytic continuation along any curve $p \rightarrow q$ with $p,q \in M$. Note that while the metrics $\psi_1^*(g)$, $\psi_2^*(g')$ are \textit{not real-analytic} in a neighborhood of the boundary (indeed, by assumption they are only log-analytic), we can consider the restrictions of these metrics to $U' = (\eta/2, \eta) \times \partial M$, which is a finite distance away from the boundary. Now because Einstein metrics are real-analytic away from the boundary (this holds in harmonic coordinates, see \cite{anderson2008survey} among other sources), we perform an analytic continuation $\chi$ along a curve starting at $p \in U' = (\eta/2, \eta)_x \times \partial M$ and ending at some other point in $M$. \nl 
\indent Let $\gamma: p \rightarrow q$ be any curve in $M$ along which we perform analytic continuation of $\chi$. We need to show that this continuation is independent of $\gamma$. Since we assume that $\pi_1(\overline{M}, \partial M) = 0$, $\gamma$ can be deformed to a curve inside of $\partial M$ and hence, inside of $U$. Since $\chi$ is determined exactly on $U \supseteq \partial M$, the extension of $\chi$ must be path-independent and we have constructed a global diffeomorphism, $\chi: \overline{M} \to \overline{M}$ such that $\chi^*(g') = g$. \nl
\indent For $m = 2m^* \geq 4$, we note that the extra assumption of $c(g) = c(g')$ replaces the determination of the conformal infinity which was previously afforded by theorem \ref{ConformalInfinityThm} in the $m = 2$ case. The rest of the proof is the same using \ref{ExpansionThm} to recover the higher $\{\omega_k, \omega_{ab}\}$ coefficients, and then constructing the same isometry as in the $m=2$ case. \qed \nl \nl
\noindent Our second set of applications uses renormalized volume in the setting of Poincar\`e-Einstein metrics. In particular, theorem \ref{ExpansionThm} determines the non-local coefficient for these metrics which gives corollary \ref{DtoNDetermination}
\begin{corollary*}
Suppose $(M^{n+1},g)$ is a Poincar\`e-Einstein manifold and $\omega_0$ is known in equation \eqref{AHExpansion}. Then knowledge of the renormalized area on all minimal surfaces for any dimension $2m^* < n + 1$ determines the non-local coefficients, $\omega_{n}$ when $n$ even and $\omega_{n-1}$ when $n$ odd, in equations \eqref{PEExpansionEven} \eqref{PEExpansionOdd}.
\end{corollary*}
\noindent \Pf Note Poincar\`e-Einstein metrics are partially even of order $2m^*$ for any $2m^* < n + 1$. Thus theorem \ref{ExpansionThm} applies directly. \qed \nl \nl
%
%
\rmk: As described in \cite{anderson2005geometric} equation 1.9, this coefficient is a major obstruction to understanding the AdS/CFT correspondence: given $(\partial M, \omega_0)$, the entire asymptotic expansion in equations \eqref{PEExpansionEven} \eqref{PEExpansionOdd} is determined by $\omega_{n-1}$ (or $\omega_n$ depending on parity). As a result, theorem \ref{ExpansionThm} tells us that features of the boundary CFT do in fact determine the non-local term. \nl \nl
Our second application lies in the case of $n = 3$, for which it is well known that $M^3 \cong \H^3 / \Gamma$ for $\Gamma$ a convex, cocompact group subgroup and the metric expands as 
\begin{equation} \label{3DPEExpansion}
g = \frac{dx^2 + \omega_0 + x^2 \omega_2 + x^4 \omega_4}{x^2}
\end{equation}
with $(\omega_4)_{ij} = (\omega_2)_{ik} (\omega_2)^k_j$ (see \cite{fefferman2011ambient}, theorem $3.4$, among other sources). In particular, the above expansion holds in an open neighborhood of the boundary $U \supseteq (\partial M) \times [0, a)$ for $a$ sufficiently small. We can also connect this to the topology of the boundary manifold in corollary \ref{3DApp}:
\begin{corollary*}
Suppose $(M^3,g)$ is a Poincar\`e-Einstein manifold with $\Sigma \cong \partial M$ and $\omega_0$ is known in equation \eqref{3DPEExpansion}. Then knowledge of the renormalized area on all minimal surfaces determines the conformal structure on $\Sigma$ and hence, $g$ globally.
\end{corollary*}
\noindent \Pf Knowledge of $\omega_0$ and $\omega_2$ determines the a conformal structure on $\Sigma$ as follows: since $g$ is known exactly in a neighborhood of the boundary via equation \eqref{3DPEExpansion}, converting the boundary defining function to $\rho$ via $\rho = \frac{1}{2} \ln(x)$ will determine the first and second fundamental form of $\Sigma_{\rho}$, which is the surface a distance $\rho$ from $\Sigma$ (see equation \cite{papadopoulos2007handbook}). These coefficients will then determine the ``first fundamental form at infinity," $I^*$ (see \cite{papadopoulos2007handbook}, Chapter 14, lemma 3.1), which is the natural conformal structure for $\Sigma$. This in turn determines $g$ (\cite{papadopoulos2007handbook}, Chapter 14, \S 1.4). \qed
%

\section{Appendix}
\subsection{Computing $L_2(\dot{g})$} \label{L2Append}
Here $L_2$ is given by the linearization if we fix a surface, $Y^{2m^*}$, and vary the metric. Then
\begin{align*}
	H^{\delta}_Y &= g^{ab}_{\delta} g_{\delta}(\nabla_a^{\delta} \p_b, \nu_{Y, \delta})
\end{align*}
Here $\{\p_a\}$ is a fixed basis for $TY$ which has vanishing christoffels at $p$ w.r.t. $g_0$. Moreover, let $\nu$ be the normal to $Y$ w.r.t. $g_0$. Then we write
\begin{align*}
	\nu_{\delta} &= A^c(\delta) \p_c + A^{\nu}(\delta) \nu \\
	A^c(0) &= 0 \\
	A^{\nu}(0) &= 1
\end{align*}
In particular $\{\p_a, \nu\}$ gives a coordinate fermi basis at $p \in Y$. Then
\begin{align*}
	g_{\delta}(\nabla_a^{\delta} \p_b, \nu_{\delta}) &= A^c(\delta) g_{\delta}(\nabla_a^{\delta} \p_b, \p_c) + A^{\nu}(\delta) g_{\delta}(\nabla_a^{\delta} \p_b, \nu) \\
	g_{\delta}(\nabla_a^{\delta} \p_b, \p_c) &= \frac{1}{2} [ \p_a g_{\delta, bc} + \p_b g_{\delta, ac} - \p_c g_{\delta, ab}] \\
	g_{\delta}(\nabla_a^{\delta} \p_b, \nu) &= \frac{1}{2} [ \p_a g_{\delta, b \nu} + \p_b g_{\delta, a \nu} - \p_{\nu} g_{\delta, ab} ] 
\end{align*}
where $\n^{\delta}$ denotes the connection on $Y$ w.r.t. $g_{\delta} \Big|_{TY}$. We differentiate
\begin{align*}
	\frac{d}{d \delta} H^{\delta}_Y \Big|_{\delta = 0} &= \frac{d}{d \delta} \left[ g^{ab}_{\delta} g_{\delta}(\nabla_a^{\delta} \p_b, \nu_{Y, \delta}) \right]_{\delta = 0} \\
	&= \dot{g}^{ab} A_{ab} + g^{ab} \frac{d}{d \delta} \left[ g_{\delta}(\nabla_a^{\delta} \p_b, \nu_{Y, \delta}) \right]_{\delta = 0} \\
	&= - \langle \dot{g}, A \rangle + g^{ab}\frac{d}{d \delta} \left[ g_{\delta}(\nabla_a^{\delta} \p_b, \nu_{Y, \delta}) \right]_{\delta = 0}  \\
	\frac{d}{d \delta} \left[ g_{\delta}(\nabla_a^{\delta} \p_b, \nu_{Y, \delta}) \right]_{\delta = 0}  &=  \dot{A}^{\nu} A_{ab} + \frac{1}{2} [ \p_a \dot{g}_{b \nu} + \p_b \dot{g}_{a \nu} - \p_{\nu} \dot{g}_{ab} ]
\end{align*}
Here, we've noted that $A^c(0) = 0$ and $g(\n_a \p_b, \p_c) = 0$ and so the first term in the expansion of	$g_{\delta}(\nabla_a^{\delta} \p_b, \nu_{\delta})$ does not contribute to the first derivative in $\delta$. To compute $\dot{A}^{\nu}$, we note that 
\begin{align*}
	g(\nu_{\delta}, \nu) &= A^{\nu}(\delta) \\
	\implies \frac{d}{d \delta} g(\nu_{\delta}, \nu) &= \dot{A}^{\nu} \\
	g_{\delta}(\nu_{\delta}, \nu_{\delta}) &= 1 \\
	\implies \dot{g}(\nu, \nu) + 2 g(\dot{\nu}, \nu) = 0 \\
	\implies \dot{A}^{\nu} &= - \frac{1}{2} \dot{g}_{\nu \nu}
\end{align*}
So that
\begin{align} \label{DerivOfMCwrtMetric}
	\frac{d}{d \delta} H_Y^{\delta} \Big|_{\delta = 0} &= -\langle \dot{g}, A \rangle - \frac{1}{2} \dot{g}_{\nu \nu} H_{Y} + \frac{1}{2} g^{ab} [\p_a \dot{g}_{b \nu} + \p_b \dot{g}_{a \nu} - \p_{\nu} \dot{g}_{ab}] \\ \nonumber
	&= - \langle \dot{g}, A \rangle + g^{ab} \p_a \dot{g}_{b \nu} - \frac{1}{2} \p_{\nu} \tr_Y(\dot{g}) \\ \nonumber
	&= - \langle \dot{g}, A \rangle + \div_Y(\dot{g}(\nu, \cdot)) - \frac{1}{2} \nu(\tr_Y(\dot{g})) \\ \nonumber
	&=: L_2(\dot{g})
\end{align}
In the first line, we recognized that $H_Y = 0$. \nl \nl
We further note that $L_2 = x^2 \overline{L}_2$ where $\overline{L}_2$ is an edge operator. To see this, first note that contracting with $\nu$ and differentiating with respect to $\nu$ are parity preserving to order $2m^*$ from equation \eqref{NormalExpansion} and the discussion in \S \ref{MinSurfBackground}. Similarly, the fact that $Y$ has an even graphical expansion up to order $2m^*$ (see equation \eqref{UExpansion}) means that computing the trace is $O(x^2)$ and parity preserving to order $2m^* + 2$. As a result 
\[
P(\dot{g}) := \div_Y(\dot{g}(\nu, \cdot)) - \frac{1}{2} \nu(\tr_Y(\dot{g}))
\]
is a first order edge operator. Now we recall that $A_{ab}$ is even to order $2m^*$ as well, which follows from computation or \cite{marx2021variations} theorem 4.1. We then have
\[
\langle \dot{g}, A \rangle = h^{ab} h^{cd} \dot{g}_{ac} A_{bd}
\]
is a $0$th order edge operator and also parity preserving to order $2m^*$. 
\subsection{Solving $J_Y(\dot{\psi}) = 0$ for $Y = HS^{m}$} \label{SolvingJacobiGeneral}
%
In this section, we explicitly find solutions to the Jacobi equation for the geodesic hemisphere in hyperbolic space via separation of variables. Given $HS^m \subseteq \H^{m+1}$, we consider variations of the boundary curve $ \gamma \mapsto \gamma_t$ with $Y_t \subseteq \H^{m+1}$ the corresponding minimal surface asymptotic to $\gamma_t$ for each $t$. We compute the Jacobi operator
\[
J_Y = \Delta_Y + ||A_Y||^2 + \Ric_g(\nu, \nu) = \Delta_Y - m
\]
when the ambient manifold is $\H^{m+1}$. We recall our parameterization for $HS^m$ from section \S \ref{HemisphereParamSection}
\begin{align*}
HS^m &= \{(x, y_1, \dots, y_m) = (x, \sqrt{1 - x^2} F_m(\theta))\} \\
v_x &= \p_x - \frac{x}{\sqrt{1 - x^2}} F_m(\theta) \\
v_{\theta_i} &= \sqrt{1 - x^2} F_{m, i}(\theta) \\
h_{xx} &= \frac{1}{(1 - x^2)x^2} \\
h_{i j} &= \frac{1 - x^2}{x^2} \bh_{i j} \\
h_{x i} &= 0 \\
g(\nabla_{v_x} v_x, v_x) &= \frac{1}{2} v_x h_{xx}  \\
&= -\frac{1 - 2x^2}{(1 - x^2)^2 x^3} \\
g(\nabla_{v_x} v_x, v_{\theta_i}) &= - \frac{1}{2} v_{\theta_i} h_{xx} \\
&= 0 \\
g(\nabla_{v_{\theta_i}} v_{\theta_j}, v_{\theta_k}) &= \frac{1}{2} [v_{\theta_i} h_{jk} + v_{\theta_j} h_{ik} - v_{\theta_k} h_{ij}] \\
&= \frac{1 - x^2}{x^2} \bGamma_{ijk}\\
g(\n_{v_{\theta_i}} v_{\theta_j}, v_x) &= - \frac{1}{2} v_x h_{ij}  \\
&= \frac{1}{x^3} \bh_{ij}
\end{align*}
where $\bh_{ij}$ and $\bGamma_{ijk}$ are the metric coefficients and Christoffel symbols on $S^m$ with respect to the euclidean metric. As a result, we have 
\begin{align*}
	\Delta_{HS^m} &= h^{\alpha \beta} [\nabla_{\alpha} \nabla_{\beta} - \nabla_{\nabla_{\alpha} \beta}] \\
	&= h^{xx} [ \p_x^2 - \nabla_{v_x} v_x] + h^{ij} [ \p_{\theta_i} \p_{\theta_j} - \nabla_{v_{\theta_i}} v_{\theta_j}]\\
	&= x^2(1 - x^2) [ \p_x^2 - h^{xx} g(\n_{v_x} v_x, v_x) \p_x] \\
	&+ \frac{x^2}{1 - x^2} \bh^{ij} [\p_{\theta_i} \p_{\theta_j} - h^{xx} g(\nabla_{v_{\theta_i}} v_{\theta_j}, v_x) \p_x - h^{k \ell} g(\n_{v_{\theta_i}} v_{\theta_j}, v_{\theta_k}) v_{\theta_{\ell}}] \\
	&= x^2(1 - x^2) \p_x^2 + x(1 - 2x^2) \p_x \\
	& + \frac{x^2}{1 - x^2} \Delta_{S^{m-1}, \bg} - h^{xx} \frac{x^2}{1 - x^2} \bh^{ij} g(\n_{v_{\theta_i} } v_{\theta_j}, v_x) \p_x \\
	& = x^2(1 - x^2) \p_x^2 + x(1 - 2x^2) \p_x + \frac{x^2}{1 - x^2} \Delta_{S^{m-1}, \bg} - x^2(1 - x^2) \frac{1}{x(1 - x^2)} \bh^{ij} \bh_{ij} \p_x \\
	& = x^2(1 - x^2) \p_x^2 + x(1 - 2x^2) \p_x + \frac{x^2}{1 - x^2} \Delta_{S^{m-1}, \bg} - x (m-1) \p_x \\ 
	&= x^2(1 - x^2) \p_x^2 + x[(2 - m) - 2x^2] \p_x + \frac{x^2}{1 - x^2} \Delta_{S^{m-1}, \bg} \\
	&= (x \p_x)^2 + (1-m) (x \p_x) + \left[ (-2x^2)x\p_x - (x^4) \p_x^2 + \frac{x^2}{1 - x^2} \Delta_{S^{m-1}, \bg} \right] \\
	\implies J_Y &= \left[ \left( (x \p_x)^2 - (m-1)(x \p_x) - m \right) + \left( (-2x^2) x\p_x - (x^4) \p_x^2  \right) \right] +  \frac{x^2}{1 - x^2} \Delta_{S^{m-1}, \bg} \\
	&= L + \frac{x^2}{1 - x^2} \Delta_{S^{m-1}, \bg}
\end{align*}
To find solutions to the Jacobi equation, we proceed with separation of variables, i.e. suppose we have a solution of the form $A(\vec{\theta}) B(x)$ where $\vec{\theta}$ parameterizes $S^{m-1}$ as in \S \ref{HemisphereParamSection}. Then
\begin{align} \nonumber
	J_Y(A B) &= 0 \\ \nonumber
	&= L(B) A + B \frac{x^2}{1 - x^2} \Delta_{S^{m-1}} \\ \label{ODEEquation}
	\implies \frac{\Delta_{S^{m-1}} A}{A} &= \frac{x^2 - 1}{x^2} \frac{L(B)}{B} = -\lambda
\end{align}
For $\lambda$ a constant. From classical techniques, we know that the eigenvalues of the laplacian on $S^{m-1}$ are negative and given (in the notation of the above) by $\lambda = \lambda_{m-1}(k) = k(m + k - 2)$. Moreover, the eigenfunctions are given by $k$th order harmonic, homogeneous polynomials on $\R^m$ restricted to $S^{m-1}$. If we parameterize 
\[
S^{m-1} = \{(\ty_1, \dots, \ty_m) = (\cos \phi, \sin \phi F_{m-2}(\theta_{m-2})) \}
\] 
and set $k = 2$ so that (see \S \ref{AngularComputation} for more details)
\begin{align} \label{SecondHarmonicEqn}
A_2(\phi, \theta_{m-2}) &= (\cos \phi)^2 - \frac{1}{m-1} (\sin \phi)^2 
\end{align}
On the RHS of equation \eqref{ODEEquation}, we have 
\begin{align} \nonumber
	\frac{x^2 - 1}{x^2} \frac{L(B)}{B} &= -\lambda_{m-1}(k) \\ \nonumber
	\implies L(B) - \frac{\lambda_{m-1}(k) x^2}{1 - x^2} B &= 0 \\ \label{BODE}
	\implies \left[ \left( (x \p_x)^2 B - (m-1)(x \p_x) B - m B \right) + \left( -x^4 \p_x^2 B + x(-2x^2)\p_x B \right) \right] - \frac{[k(m+k-2)] x^2}{1 - x^2} B & = 0 
\end{align}
\subsubsection{Solving when $k = 2$ on $HS^{2m^*}$} \label{Solvingk2}
We consider the case of set $k = 2$ and $m = 2m^*$. If we write $B = B_{2,2m^*} = \frac{\overline{B}_{2,2m^*}}{x}$, then
\begin{align*}
	[(x\p_x)^2 - (m-1)(x\p_x) - m](B_{2,2m^*}) &= x(\overline{B}_{2,2m^*})_{xx} - m(\overline{B}_{2,2m^*})_{x} \\
	\left( -x^4 \p_x^2  + x(-2x^2)\p_x \right) B_{2,2m^*} &= -x^3 (\overline{B}_{2,2m^*})_{xx} 
\end{align*}
We now impose $\overline{B}_{2,2m^*}(0) = 1$ to get
\begin{align} \label{B2ODE}
L_{2,2m^*}(\bB_{2,2m^*}) &= (x - x^3) (\overline{B}_{2,2m^*})_{xx} - 2m^* (\overline{B}_{2,2m^*})_{x} - \frac{4m^* x}{1 - x^2} \overline{B}_{2,2m^*} = 0
\end{align}
We note that the indicial roots of this equation at $x = 0$ are given by $r = 0, 2m^* + 1$, and at $x = 1$ are given by $r = 1, -m^*$. Writing $\overline{B}_{2,2m^*} = \frac{f(x)}{(1-x^2)^{m^*}}$, this produces a hypergeometric differential equation for which the general solution is
\begin{align} \label{B2Equation}
\overline{B}_{2,2m^*}(x) &= \frac{c_1 \cdot \F21(-m^*, -\frac{1}{2}-m^*, \frac{1}{2}-m^*, x^2) + c_2 \cdot x^{2m+1}}{(1-x^2)^{m^*}}
\end{align}
where $\F21$ is the hypergeometric function. We record a few important properties of $\F21$:
\begin{align*}
F_1(a,b,c,0) &= 1 \\
F_1(a,b,c,1) &= \frac{\Gamma(c) \Gamma(c-a-b)}{\Gamma(c-a) \Gamma(c-b)}  \quad \text{when} \quad \text{Re}(c-b-a) > 0
\end{align*}
where $\Gamma(n+1) = n!$ for $n \in \Z^{\geq 0}$ and in general $\Gamma(z+1) = z \Gamma(z)$. Note that to enforce $\overline{B}_2(0) = 1$ we have that $c_1 = 1$. To have $\overline{B}_2$ be well defined on all of $[0,1]$, we are forced to set
\begin{equation} \label{c2Computation}
c_2:= - F_1(-m^*, -\frac{1}{2} - m^*, \frac{1}{2} - m^*, 1) = - \frac{\Gamma(1+m^*)\Gamma(1/2 - m^*)}{\Gamma(1/2) \Gamma(1)} = - \frac{(m^*)!}{\left(-\frac{1}{2} \right) \cdots \left(\frac{1}{2} - m^* \right)}
\end{equation}
While this choice of $c_2$ only makes the numerator of equation \eqref{B2Equation} vanish at $x = 1$ a priori, it follows from the Fuchsian theory of ODEs with rational coefficients (see \cite{gray1984fuchs} among many sources) that this choice of $c_2$ is the unique solution of equation \eqref{B2ODE} such that $B_2(0) = 1$ and $B_2(1) = 0$ in accordance with our analysis of the indicial roots of $r = 0$ at $x = 0$ and $r = 1$ at $x = 1$. \nl \nl
As an example, when $m^* = 1$ and $k = 1$, we have
\begin{align*}
	(x - x^3) (\overline{B}_{2,2})_{xx} - 2 (\overline{B}_{2,2})_x - \frac{4x}{1 - x^2} \overline{B}_{2,2} &= 0  \\
	\implies \overline{B}_{2,2}(x) = c_1 \frac{1 - 3x^2}{x^2 - 1} + c_2\frac{x^3}{x^2 - 1}
\end{align*}
again, enforcing $\overline{B}_{2,2}(0) = 1$ and $\overline{B}_{2,2}(1)$ to be defined, we see that $c_1 = 1$ and $c_2 = 2$ (agreeing with equation \eqref{c2Computation}) so
\[
\overline{B}_{2,2}(x) = \frac{(1 - x)(2x + 1)}{1 + x} 
\]
which is defined everywhere on $[0,1]$. The corresponding Jacobi field is then 
\begin{align*}
\dot{\psi}_2(\phi, x) &= \cos(2 \phi) \frac{\overline{B}_{2,2}(x)}{x} \\
&= \cos(2 \phi) \frac{(1 - x)(2x + 1)}{x(1 + x)}
\end{align*}
\subsection{Computation of finite part} \label{FinitePartComputation}
The goal of this section is to show that the following quantity
\begin{equation} \label{IkmIntegral}
I_{k,m^*}: = \FPz \int_0^1 x^{z+2k-2m^*}[{}_2F_1(-m-1/2,-m;1/2-m,x^2) + c x^{2m+1}]
\end{equation}
is non-zero, for $c = - {}_2F_1(-m-1/2,-m;1/2-m,1)$. To verify $R \neq 0$, we use that ${}_2F_1(-m-1/2,-m;1/2-m,x^2)$ has a finite expansion and recognize the finite part of the integral as the evaluation of \textit{another} hypergeometric function.
\begin{align*}
\int_0^1 x^{z+2k-2m^*} {}_2F_1(-m-1/2,-m;1/2-m,x^2) &= \sum_{\ell = 0}^{m} \frac{(-m-1/2)_{\ell} (-m)_{\ell}}{(1/2-m)_{\ell}} \FPz \int_0^1 \frac{x^{z+2k-2m^*+2\ell}}{\ell!} \\
&= \sum_{\ell = 0}^{m} \frac{(-m-1/2)_{\ell} (-m)_{\ell}}{(1/2-m)_{\ell}} \frac{1}{2k-2m^*+2\ell + 1}  \frac{1}{\ell!}
\end{align*}
The key is to now view this as a hypergeometric function of its own by writing
\begin{align*}
\frac{1}{2k-2m^*+2\ell + 1}  &= \frac{1}{2[k-m^*+ (\ell - 1) + 3/2]}  \\
&= \frac{(k-m^* + 1/2)_{\ell}}{(k-m^* + 3/2)_{\ell}} \\
\sum_{\ell = 0}^{m} \frac{(-m-1/2)_{\ell} (-m)_{\ell}}{(1/2-m)_{\ell}} \frac{1}{2k-2m^*+2\ell + 1}  \frac{1}{\ell!} &= \sum_{\ell = 0}^{m} \frac{(-m-1/2)_{\ell} (-m)_{\ell}}{(1/2-m)_{\ell}} \frac{(k-m^* + 1/2)_{\ell}}{(k-m^* + 3/2)_{\ell}}  \frac{1}{\ell!}  \\
& = {}_3F_2\left( -m-1/2,-m, k - m^* + 1/2; 1/2-m, k-m^* + 3/2; 1\right)
\end{align*}
where $\F32$ is a generalized hypergeometric function. We can compute this as follows:
\begin{align*}
\frac{(-m^*-1/2)_{\ell} }{(1/2-m^*)_{\ell}} \frac{(k-m^* + 1/2)_{\ell}}{(k-m^* + 3/2)_{\ell}} &= \frac{(-m^*-1/2)(k-m^* + 1/2)}{(-m^*-1/2 + \ell) (k-m^* + 1/2 + \ell)} \\
&= \frac{(-m^*-1/2)(k-m^* + 1/2)}{k+1} \left[\frac{1}{-m^* - 1/2 + \ell} - \frac{1}{k - m^* + 1/2 + \ell} \right] \\
&= \frac{1}{k+1} \left[\frac{(-m^* - 1/2)_{\ell} (k-m^* + 1/2)}{(1/2-m^*)_{\ell}} - \frac{(-m^*-1/2)(k-m^* + 1/2)_{\ell}}{(k-m^* + 3/2)_{\ell}} \right]
\end{align*}
which means that we can break our ${}_3F_2$ expression into the following two ${}_2 F_1$ expressions:
\begin{align*}
{}_3F_2\left( -m-1/2,-m, k - m^* + 1/2; 1/2-m, k-m^* + 3/2; 1\right) & = \frac{k-m^* + 1/2}{k+1} {}_2F_1(-m^* - 1/2, -m^*; 1/2 - m^*; 1) \\
& - \frac{(-m^*-1/2)}{k+1}{}_2 F_1(k-m^* + 1/2, -m^*; k - m^* + 3/2; 1)
\end{align*}
we can compute these explicitly via combinatorics, or we can use the Gauss transformation
\[
F(a,b,c,1) = \frac{\Gamma(c) \Gamma(c-b-a)}{\Gamma(c-b) \Gamma(c-a)}
\]
which means that 
\begin{align*}
{}_3F_2\left( -m-1/2,-m, k - m^* + 1/2; 1/2-m, k-m^* + 3/2; 1\right) & = \frac{k-m^* + 1/2}{k+1} \frac{\Gamma(1/2-m^*) \Gamma(m^* + 1)}{\Gamma(1) \Gamma(1/2)} \\
& - \frac{(-m^*-1/2)}{k+1}\frac{\Gamma(k-m^* + 3/2) \Gamma(m^* + 1)}{\Gamma(1) \Gamma(k + 3/2)}
\end{align*}
%
Noting that 
\[
c = - {}_2F_1(-m-1/2,-m;1/2-m,1) = -\frac{\Gamma(1/2-m^*) \Gamma(m^* + 1)}{\Gamma(1) \Gamma(1/2)}
\]
we have that 
\begin{align*}
\FPz \int_0^1 x^{z+2k-2m^*} c x^{2m^* + 1} &= \int_0^1 x^{z+2k + 1} c \\
&= \frac{c}{2(k+1)}
\end{align*}
We can now evaluate our original integral \eqref{IkmIntegral}:
\begin{align*}
I_{k,m^*} &= \FPz \int_0^1 x^{z+2k-2m^*}[{}_2F_1(-m-1/2,-m;1/2-m,x^2) + c x^{2m+1}] \\
&= \frac{k-m^* + 1/2}{k+1} \frac{\Gamma(1/2-m^*) \Gamma(m^* + 1)}{\Gamma(1) \Gamma(1/2)} - \frac{(-m^*-1/2)}{k+1}\frac{\Gamma(k-m^* + 3/2) \Gamma(m^* + 1)}{\Gamma(1) \Gamma(k + 3/2)} + \frac{c}{2(k+1)}
\end{align*}
our goal is to show that $I_{k,m^*} \neq 0$. We group like terms on the right hand side
\begin{align*}
I_{k,m^*} &= \frac{\Gamma(m^* + 1)}{k+1} \left[\frac{\Gamma(1/2-m^*)}{\Gamma(1/2)} \left( - \frac{1}{2} + (k-m^* + 1/2) \right) - \frac{(-m^* - 1/2) \Gamma(k-m^* + 3/2)}{\Gamma(k+3/2)}\right] \\
\implies \frac{2I_{k,m^*}(k+1)}{\Gamma(m^*+1)} &= \left[\frac{\Gamma(1/2-m^*)}{\Gamma(1/2)} 2(k-m^*) - \frac{(-2m^* -1) \Gamma(k-m^* + 3/2)}{\Gamma(k+3/2)}\right]
\end{align*}
To see that the right hand side is non-zero, note that 
\begin{align*}
\frac{\Gamma(1/2-m^*)}{\Gamma(1/2)} &= \frac{2^{m^*}}{(-1)(-3)\cdots (-2m^* + 1)}\\
\frac{\Gamma(k-m^* + 3/2)}{\Gamma(k+3/2)} &= \frac{2^{m^*}}{(k+1/2)(k-1/2) \cdots (k-m^*+3/2)}
\end{align*}
In particular, if we let $k-m^* = 2^r L$ where $L$ is odd and $r \geq 0$, then
\begin{align*}
\frac{\Gamma(1/2-m^*)}{\Gamma(1/2)} 2(k-m^*) & = 2^{m^* + 1 + r} Q_1 \\
\frac{(-2m^* -1) \Gamma(k-m^* + 3/2)}{\Gamma(k+3/2)} &= 2^{m^*} Q_2
\end{align*}
where $Q_1,Q_2$ are both rational numbers with odd numerators and denominators. Thus by parity considerations
\[
\frac{\Gamma(1/2-m^*)}{\Gamma(1/2)} 2(k-m^*) \neq \frac{(-2m^* -1) \Gamma(k-m^* + 1/2)}{\Gamma(k+3/2)}
\]
meaning that $I_{k,m^*}$ is non-zero!
\subsection{Angular Computation} \label{AngularComputation}
The variations to be used in equation \eqref{DifferentiateEqn} are given by:
\[
\dot{\phi}(\vec{\theta}, x) = \frac{1}{x} \bB_{2,2m^*}(x) f_{2,2m^*-1}(\vec{\theta})
\]
where 
\begin{align*}
\bB_{2,2m^*}'' - \frac{2m^*}{x(1-x^2)} \bB_{2,2m^*}' - \frac{4m}{(1-x^2)^2} \bB_{2,2m^*} &= 0 \\
\bB_{2,2m^*}(x) &= \frac{{}_2F_1(-m^*, -m^* - 1/2; 1/2-m^*; x^2) + c x^{2m^* + 1}}{(1-x^2)^{m^*}} \\
\Delta_{S^{2m^*-1}} f_{2, 2m^*-1}(\vec{\theta}) &= -4m^*  f_{2,2m^*-1}
\end{align*}
Such eigenfunctions are restrictions of harmonic, homogeneous polynomials on $\R^{2m^*}$, restricted to $S^{2m^*-1}$. In this case, we parameterize 
\begin{align*}
S^{2m^* - 1} &= \{(y_1, y_2, \dots, y_{2m}) = (\cos \phi, \sin \phi F_{2m^*-2}(\theta_{2m^*-2})) \} \subseteq \R^{2m^*}
\end{align*}
where $F_{2m^*-2}(\theta_{2m^*-2})$ is a polar coordinate parameterization of $S^{2m^*-2} \subseteq \R^{2m^*-1}$. In these coordinates
\[
f_2(\phi, \theta_{m-2}) = (\cos^2 \phi) - \frac{1}{2m^*-1} (\sin^2 \phi) = \frac{y_1^2}{r^2} - \frac{1}{(2m^*-1)r^2}\left[y_2^2 + \dots y_{2m^*}^2 \right]
\]
The goal of this section is to evaluate the angular integrals arising in the evaluation of $\frac{d}{dt} L(t)$ from equation \eqref{DifferentiateEqn} from section \S \ref{GeneralMProof}. The two integrals of interest are given by \eqref{I1Integral} \eqref{I2Integral}:
\begin{align*}
I_1 &= \int_{S^{2m^*-1}} \left[ \cos^2(\phi) - \frac{1}{m-1} \sin^2(\phi) \right] \omega_k(F_{2m^*-1}(\theta_{2m^*-1}), F_{m-1}(\theta_{2m^*-1})) dA_{S^{2m^*-1}}(\theta_{2m^*-1}) \\
I_2 &= \int_{S^{2m^*-1}} \left[ \cos^2(\phi) - \frac{1}{2m^*-1} \sin^2(\phi) \right] (\tr_{S^{2m^*-1}}\omega_k) (\sin^{2m^*-2} \phi)  d \phi dA_{S^{2m^*-2}}(\theta_{2m^*-2}) 
\end{align*}
We can briefly reduce
\begin{align*}
I_1 &= \int_{S^{2m^*-1}} \left[ \cos^2(\phi) - \frac{1}{2m^*-1} \sin^2(\phi) \right] \omega_k(F_{2m^*-1}(\theta_{2m^*-1}), F_{2m^*-1}(\theta_{2m^*-1})) dA_{S^{2m^*-1}}(\theta_{2m^*-1}) \\
&= \int_{\phi = 0}^{\pi} \left[ \cos^2(\phi) - \frac{1}{2m^*-1} \sin^2(\phi) \right] \omega_k((\cos \phi, \sin \phi F_{2m^*-2}(\theta_{2m^*-2})), (\cos \phi, \sin \phi F_{2m^*-2}(\theta_{2m^*-2}))) \sin^{2m^*-2} \phi d \phi dA_{S^{2m^*-2}} \\
&= \int_{\phi = 0}^{\pi} \left[ \cos^2(\phi) - \frac{1}{2m^*-1} \sin^2(\phi) \right] \left[ (\cos^2 \phi)\omega_{k,11} + 2 \sin\phi \cos \phi \omega_{k}(\p_{y_1}, F_{2m^*-2}(\theta)) + \sin^2\phi \tr_{S^{2m^*-2}}(\omega_k) \right] \sin^{2m^*-2} \phi d \phi dA_{S^{2m^*-2}}(\theta_{2m^*-2})
\end{align*}
Similarly,
\begin{align*}
I_2 &= \int_{S^{2m^*-1}} \left[ \cos^2(\phi) - \frac{1}{2m^*-1} \sin^2(\phi) \right] (\tr_{S^{2m^*-1}}\omega_k) (\sin^{2m^*-2} \phi)  d \phi dA_{S^{2m^*-2}}(\theta_{2m^*-2}) \\
&= \int_{S^{2m^*-1}} \left[ \cos^2(\phi) - \frac{1}{2m^*-1} \sin^2(\phi) \right] [(\sin^2 \phi) \omega_{k,11} + (\cos^2 \phi) \tr_{S^{2m^*-2}}(\omega_k) ] (\sin^{2m^*-2} \phi)  d \phi dA_{S^{2m^*-2}}(\theta_{2m^*-2})
\end{align*}
We now compute the $\phi$-dependent integrals
\begin{align*}
\int_{0}^{\pi} f_2(\phi) \cos^2(\phi) \sin^{2m^*-2}(\phi) &= \int_{0}^{\pi} \cos^4 \sin^{2m^*-2} - \frac{1}{2m^*-1} \cos^2 \sin^{2m^*} \\
&= (2/(2m^*-1)) \int_0^{\pi} \cos^2(\phi) \sin^{2m^*}(\phi)\\
&=: c_{2m^*-1} > 0 \\
\int_{0}^{\pi} f_2(\phi) \sin(\phi) \cos(\phi) \sin^{2m^*-2}(\phi) &= 0 \\
\int_{0}^{\pi} f_2(\phi) \sin^2(\phi) \sin^{2m^*-2} &= \int_{0}^{\pi} \cos^2 \sin^{2m^*} - \frac{1}{2m^*-1} \sin^{2m^*+2} \\
&= - c_{2m^*-1} 
\end{align*}
where we have used integration by parts on $\cos^4 \sin^{2m^*-2} = (\cos^3) (\sin^{2m^*-2} \cos)$ and $\sin^{2m^*+2} = (\sin^{2m^*+1})(\sin)$. Finally, we note that 
\[
\int_{S^{2m^*-2}} \tr_{S^{2m^*-2}}(\omega_k) dA_{S^{2m^*-2}} = \rho_{2m^*-2} [ \omega^k_{22} + \omega^k_{33} + \dots + \omega^k_{(2m^*)(2m^*)}]
\]
for $\rho_{2m^*-2} > 0$. Finally, let $\tau_{2m^*}$ denote the volume of the sphere of unit radius and dimension $2m^*$. With this, we see that 
\begin{align*}
I_1 & = \omega^k_{11} c_{2m^*-1} \tau_{2m^*-2} - c_{2m^*-1} \rho_{2m^*-2} [ \omega^k_{22} + \dots + \omega^k_{(2m^*)(2m^*)}]\\
I_2 & = -c_{2m^*-1} \omega^k_{11} \tau_{2m^*-2} + c_{2m^*-1} \rho_{2m^*-2} [ \omega^k_{22} + \dots + \omega^k_{(2m^*)(2m^*)}]
\end{align*}
And so, in particular
\[
x^{z + 2k - 2m^*} I_1 + x^{z + 2k - 2m^* + 2} I_2 = c_{2m^*-1} \left(\omega^k_{11} \tau_{2m^*-2} - \rho_{2m^*-2}[\omega^k_{22} + \dots + \omega^k_{(2m^*)(2m^*)}] \right) (1-x^2) x^{z + 2k - 2m^*}
\]
\bibliography{Renormalized_Volume_Inverse}{}
\bibliographystyle{plain}
\end{document}